\documentclass[a4paper,10pt]{article}
\usepackage{fullpage}
\usepackage{amsmath,amssymb,latexsym,amsthm}
\usepackage[pdftex]{graphicx}
\usepackage{graphics}
\usepackage{psfrag}
\usepackage{color, epsfig}
\usepackage{mathrsfs}
\usepackage{bbm}
\usepackage[colorlinks=true]{hyperref} 
\usepackage[frenchb,english]{babel}
\usepackage{tikz}
\usetikzlibrary{hobby}
\usepackage{stmaryrd}

\def\R{{\mathbb R}}
\def\N{{\mathbb N}}

\newcommand{\Hreg}{\mathbf{H}_{\omega,\Theta}}

\newcommand{\norm}[1]{\left\Vert#1\right\Vert}

\renewcommand{\leq}{\leqslant}
\renewcommand{\geq}{\geqslant}

\newtheorem{theorem}{Theorem}[section]
\newtheorem{definition}[theorem]{Definition}

\newtheorem{proposition}[theorem]{Proposition}

\newtheorem{remark}[theorem]{Remark}
\numberwithin{equation}{section}

\newcommand{\ma}[1]{{\color{black}#1}}
\newcommand{\red}[1]{{\color{black}#1}}
\title{\ma{Insensitizing control for the heat equation with respect to domain variations}\footnote{The first author has been supported by the Agence Nationale de la Recherche, Project IFSMACS, grant ANR-15-CE40-0010, and by the CIMI Labex, Toulouse, France, under grant ANR-11-LABX-0040-CIMI. The three authors are supported by the Agence Nationale de la Recherche, Project TRECOS. the third authors has been partially supported by the ANR Project ``SHAPe Optimization - SHAPO'' and by the Project ``Analysis and simulation of optimal shapes - application to lifesciences'' of the Paris City Hall. 
}}

\author{
Sylvain Ervedoza\footnote{Institut de Math\'ematiques de Bordeaux; UMR 5251; Universit\'e de Bordeaux ; CNRS ; Bordeaux INP; F-33400 Talence, France. {\tt sylvain.ervedoza@math.u-bordeaux.fr}}
\and
Pierre Lissy\footnote{Ceremade,   Universit\'e  Paris-Dauphine \&  CNRS  UMR  7534,  Universit\'e PSL,  75016  Paris,   France (\texttt{lissy@ceremade.dauphine.fr}).}
\and
Yannick Privat\footnote{Universit\'e de Strasbourg, CNRS UMR 7501, INRIA, Institut de Recherche Math\'ematique Avanc\'ee (IRMA), 7 rue Ren\'e Descartes, 67084 Strasbourg, France (\texttt{yannick.privat@unistra.fr}).}
}

\date\today
\begin{document}
\maketitle
\begin{abstract}
This article is dedicated to insensitizing issues for a quadratic functional involving the solution of the linear heat equation with respect to domain variations. This work can be seen as a continuation of [P. Lissy, Y. Privat, and Y. Simpor\'e. Insensitizing control for linear and semi-linear heat equations with partially unknown domain. ESAIM Control Optim. Calc. Var., 25:Art. 50, 21, 2019], insofar as we generalize several of the results it contains and investigate new related properties. In our framework, we consider boundary variations of the spatial domain on which the solution of the PDE is defined at each time, and investigate three main issues: (i) approximate insensitizing, (ii) approximate insensitizing combined with an exact insensitizing for a finite-dimensional subspace, and (iii) exact insensitizing. We provide  positive answers to questions (i) and (ii) and partial results to question (iii).
\end{abstract}
\medskip


\tableofcontents

\section{Introduction}

\subsection{Insensitizing controls with respect to domain variations, framework}
The goal of this article is to discuss an insensitizing control problem for the heat equation with respect to variations of the boundary. First results in this direction have already been obtained in \cite{Lissy-Privat-Simpore}. Introducing this problem precisely requires some notations, which we choose similar to the ones in \cite{Lissy-Privat-Simpore}. 

Let $T>0$ denote a horizon of time, $\omega$ and $\Theta$ be two open subsets of $\R^d$, $d\in\mathbb N^*$, and $\xi \in L^2(0,T; L^2(\R^d))$.  

For $\Omega$ a bounded, connected and open set of $\mathbb{R}^d$ of class $\mathscr{C}^2$, we consider the shape functional $J_h$, indexed by $h \in L^2(0,T; L^2( \omega))$, defined by 
\begin{equation}
	\label{def-J}
	J_h(\Omega) = \dfrac{1}{2}\int_{0}^{T} \int_{\Theta}y_{\Omega,h}(t,x)^{2}\, dxdt, 
\end{equation}
where $y_{\Omega,h}\in L^2((0,T)\times \mathbb R^d)$ is defined on $(0,T)\times\Omega$ as the unique weak solution of 
\begin{equation}
	\label{y-Omega-h}
		\left\{ 
			\begin{array}{rcll}
				\dfrac{\partial y}{\partial t}-\Delta y &=&\xi +h\mathbbm{1}_{\omega}& \text{ in }(0,T) \times \Omega,
				\\ 
				y&=& 0 & \text{ on }(0,T) \times \partial \Omega,
				\\ 
				y\left( 0,\cdot\right) &=& 0 & \text{ in }\Omega,
			\end{array}
		\right.
\end{equation}
($\mathbbm{1}_\omega$ denotes the characteristic function of the set $\omega$), extended by $0$ outside $(0,T) \times \Omega$.

Originally, the insensitizing problem consists in finding a control function such that some functional depending on the solution of a partial differential equation is locally insensitive to perturbations of the initial condition. This issue was first raised by J.L. Lions in \cite{Lions-1990}. We refer to Section~\ref{sec:SoA} for bibliographical comments. \ma{Let us mention that in the references \cite{Tebou1,Tebou2}, this question is called desensitizing problem.}
Nevertheless, up to our knowledge, insensitizing properties with respect to shape variation issues have been first investigated in \cite{Lissy-Privat-Simpore}. Let us recall here what we are talking about: given $\Omega_0$ a bounded, connected and open set of $\R^d$ with $\mathscr{C}^2$ boundary, our goal is to determine, whenever it exists, a control function $h \in L^2(0,T;L^2(\omega))$ such that $J_h$ does not depend on small variations of $\Omega$ in a neighbourhood of $\Omega_0$ (which will be made precise in what follows) at first order. In other words, we want to choose the control function $h$ in such a way that the functional $J_h$ is insensitized with respect to small variations of the domain. 

To give a precise meaning to this, we first remark that this problem is interesting only  if the intersection of these two last sets with $\Omega_0$ is nonempty, in which case the functional $J_h$ only depends on $\Omega\cap\omega$ and $\Omega\cap\Theta$. Hence, in the following, we will assume that $\Omega_0$ is a bounded, connected and open set of $\R^d$ with a $\mathscr{C}^2$ boundary and that $\omega$ and $\Theta$ are open subsets of $\Omega_{0}$. 

It is convenient to endow the set of domains with some differential structure. In what follows, we will use the notion of differentiation in the sense of Hadamard \cite{Delfour-Zolesio,Henrot-Pierre-2005}, which is classically used in the framework of shape optimization. This means that perturbations of $\Omega_0$ will be defined with the help of well-chosen diffeomorphisms, which have the advantage of preserving some topological features such as connectedness, boundedness and regularity.

Accordingly, we introduce the class $W^{3, \infty}(\R^d, \R^d)$ of admissible vector fields. Then,  for each element $\mathbf{V}$ of $W^{3,\infty}(\R^d, \R^d)$, there exists $\tau_{\mathbf{V}}>0$ such that for all $\tau\in [0,\tau_{\mathbf{V}})$, the mapping $\mathbf{T}_\tau:=\operatorname{Id}+\tau\mathbf{V}$ defines a diffeomorphism\footnote{To be more precise, it is easy to see that the choice $\tau_{\mathbf{V}}=1/\Vert \mathbf{V}\Vert_{W^{3,\infty}}$ works, by applying the Banach fixed-point theorem.} in $\R^d$, \textit{i.e.} the mapping $\mathbf{T}_\tau$ is invertible and $\mathbf{T}_\tau^{-1}\in W^{3,\infty}(\R^d,\R^d)$.
Furthermore, since $\mathbf{T}_\tau$ writes as a perturbation of the identity operator, one easily infers that $\mathbf{T_\tau}(\Omega_0)$ is a connected, bounded domain whose boundary $\mathbf{T_\tau}(\partial \Omega_0)$ is of class $\mathscr{C}^2$. It is notable that, in this framework, one has $\partial \mathbf{T_\tau}(\Omega_0)=\mathbf{T_\tau}(\partial \Omega_0)$. \ma{To sum-up, the assumption that $\mathbf{V}$ belongs to $W^{3,\infty}(\R^d, \R^d)$ is essentially technical and preserves the $\mathscr{C}^2$ regularity of the domain boundary once the domain deformation is applied \cite[Chapter~7]{Delfour-Zolesio}.}

In the sequel, given $\mathbf{V}\in W^{3,\infty}(\R^d, \R^d)$, we introduce the family $\{\Omega_{\tau\mathbf{V}}\}_{\tau\in [0,\tau_{\mathbf{V}})}$ of domains defined by
\[
	\Omega_{\tau\mathbf{V}}= (\operatorname{Id}+\tau\mathbf{V})(\Omega_0).
\]
As a consequence of the above discussion, for $\tau \in [0,\tau_{\mathbf{V}})$, each domain $\Omega_{\tau\mathbf{V}}$ inherits the aforementioned properties.

It is then classical (see \textit{e.g.} \cite[Chap. 5]{Henrot-Pierre-2005}) that the map $\tau \mapsto J_h(\Omega_{\tau\mathbf{V}})$ is differentiable in a neighbourhood of $\tau = 0$. In the following result, we provide a workable expression of the shape derivative $\left.dJ_h(\Omega_{\tau\mathbf{V}})/d\tau\right|_{\tau=0}$.
\begin{proposition}[{\cite[Proof of Proposition 1.1]{Lissy-Privat-Simpore}}]
	\label{Prop-J-differentiable}
	Let $\xi \in L^2(0,T; L^2(\R^d))$ and $h \in L^2(0,T; L^2(\omega))$. For all $ \mathbf{V} \in W^{3, \infty} (\R^d, \R^d)$, the mapping $\tau \mapsto J_h(\Omega_{\tau\mathbf{V}})$ is differentiable at $\tau = 0$ and 
	\begin{align}
%		\label{Diff-J-V}
		\left.\dfrac{d }{d \tau}\left({J}_h\left(\Omega_{\tau\mathbf{V}} \right)\right)\right|_{\tau=0}
		%&= 		\int_0^T \int_\Theta y_0 (t,x) y_{\mathbf{V}}(t,x) \, dx dt\\
		\label{Diff-J-V-2}
		& = 
		\int_{\partial\Omega_0} \mathbf{V} \cdot \mathbf{n} \left(  \int_0^T \partial_n y_0 \partial_n z_0 \, dt\right) d\sigma , 
	\end{align}
	where the pair $(y_0,z_0)$ solves the coupled forward-backward system 
	\begin{equation}
		\label{y-Omega-0-h}
		\left\{ 
			\begin{array}{rcll}
				\dfrac{\partial y_0}{\partial t}-\Delta y_0 &=&\xi +h\mathbbm{1}_{\omega
}& \text{ in }(0,T) \times \Omega_0,
				\\ 
				y_0 &=& 0 & \text{ on }(0,T) \times \partial \Omega_0,
				\\ 
				y_0\left( 0,\cdot\right) &=& 0 & \text{ in }\Omega_0,
\\
				-\dfrac{\partial z_0}{\partial t}-\Delta z_0 &=&\mathbbm{1}_{\Theta}y_0 & \text{ in }(0,T) \times \Omega_0,
				\\ 
				z_0 &=& 0 & \text{ on }(0,T) \times \partial \Omega_0,
				\\ 
				z_0\left( T,\cdot\right) &=& 0 & \text{ in }\Omega_0.
			\end{array}
		\right.
	\end{equation}
	($\mathbbm{1}_\Theta$ denotes the characteristic function of the set $\Theta$.)
\end{proposition}

We now recall the precise definitions of insensitizing that will be used next, introduced in \cite[Definition 1.1]{Lissy-Privat-Simpore}\footnote{In this reference, the authors restricted \red{the properties below} to diffeomorphisms $\mathbf{V} \in  W^{3,\infty}(\R^d, \R^d)$ of norm less than $1$, but an easy homogeneity argument enables to give an equivalent definition for any $\mathbf V  \in  W^{3,\infty}(\R^d, \R^d)$.} and much inspired of notions introduced in \cite{Lions-Sentinelles-1992,Bodart-Fabre}.
\begin{definition}
Let {$\xi \in L^2(0,T; L^2(\R^d))$}.
	%\label{def:Insensitizing}
	%
	We say that the control function $h \in L^2(0,T; L^2(\omega))$ insensitizes $J_h$ exactly at $\Omega_0$ at the first order with respect to boundary perturbation if
	\begin{equation}
		\label{Exact-Insentitizing}
		\text{for all } \mathbf{V}\in W^{3,\infty}(\R^d, \R^d),
		\quad
		\left.\dfrac{d }{d \tau}\left(J_h \left(\Omega_{\tau\mathbf{V}} \right)\right)\right|_{\tau=0}=0.
	\end{equation}

	Let $\mathcal{E}$ be a linear subspace of $W^{3,\infty}(\R^d, \R^d)$. 
	We say that the control function $h \in L^2(0,T; L^2(\omega))$ exactly insensitizes $J_h$ for  $\mathcal{E}$ at $\Omega_0$ at the first order with respect to boundary perturbation if
	\begin{equation}
		\label{Finite-Insensitizing}
		\text{for all } \mathbf{V}\in \mathcal{E},
		\quad
		\left.\dfrac{d }{d \tau}\left(J_h \left(\Omega_{\tau\mathbf{V}} \right)\right)\right|_{\tau=0}=0.
	\end{equation}

	Given $\varepsilon>0$, we say that the control function $h \in L^2(0,T; L^2(\omega))$ $\varepsilon$-approximately insensitizes $J_h$ at $\Omega_0$ at the first order with respect to boundary perturbation if
	\begin{equation}
		\label{Approx-Insensitizing}
		\text{for all }  \mathbf{V}\in W^{3,\infty}(\R^d, \R^d),
		\quad
		 \left|\left.\dfrac{d }{d \tau}\left(J_h\left(\Omega_{\tau \mathbf{V}}\right)\right)\right|_{\tau=0}\right |\leq \varepsilon \| \mathbf{V} \|_{W^{3, \infty}(\R^d; \R^d)}.
\end{equation}
\end{definition}

Let us conclude this section by introducing interesting issues related to insensitizing of the solution of the heat equation with respect to domain variations, that will be tackled in what follows:
\textsf{\begin{enumerate}
\item[\bf Q1.]  ({\bf $\varepsilon$-approximate insensitizing}) Let $\xi \in L^2(0,T; L^2(\R^d))$ and $\varepsilon>0$. Does there exist a control function $h \in L^2(0,T; L^2(\omega))$ that $\varepsilon$-approximately insensitizes $J_h$ at $\Omega_0$ in the sense of \eqref{Approx-Insensitizing}?
\item[\bf Q1'.]  ({\bf $\varepsilon$-approximate insensitizing and null/approximate controllability}) If the answer to {\bf Q1} is yes, is it possible to choose  $h \in L^2(0,T; L^2(\omega))$ in such a way that it is also {a} null control\footnote{\red{The wording ``null control'' refers to a function $h$ such that the solution $(y_0,z_0)$ of \eqref{y-Omega-0-h} satisfies $y_0(T) = 0$ in $\Omega_0$. }} or a  $\varepsilon$-approximate control\footnote{\red{Given $y_T\in L^2(\Omega_0)$ and $\varepsilon>0$, the wording ``$\varepsilon$-approximate control'' refers to a control function $h \in L^2(0,T; L^2(\omega))$ such that the solution $(y_0,z_0)$ of \eqref{y-Omega-0-h} satisfies $\norm{y_0(T) - y_T}_{L^2(\Omega_0)} \leq \varepsilon$.}} for $y_0$ at time $T$?
\item[\bf Q2.] ({\bf $\varepsilon$-approximate insensitizing and exact insensitizing for a finite-dimensional subspace \red{$\mathcal{E}$}})  Let {$\xi \in L^2(0,T; L^2(\R^d))$} and $\varepsilon>0$. Does there exist a control $h\in L^2(0,T; L^2(\omega))$ that insensitizes $J_h$ exactly for $\mathcal{E}$ in the sense of \eqref{Finite-Insensitizing}, and at the same time that $\varepsilon$-approximately insensitizes $J_h$ in the sense of \eqref{Approx-Insensitizing}?
 \item[\bf Q2'.]  ({\bf $\varepsilon${-}approximate insensitizing, exact insensitizing for a finite-dimensional subspace \red{$\mathcal{E}$}, and null/approximate controllability}) If the answer to {\bf Q2} is yes, is it possible to choose  $h \in L^2(0,T; L^2(\omega))$ in such a way that it is also {a} null control or a $\varepsilon$-approximate control for $y_0$ at time $T$, in the sense given in {\bf Q1'}?
\item[\bf Q3.] ({\bf exact insensitizing})  Let {$\xi \in L^2(0,T; L^2(\R^d))$}. Is it possible to exactly insensitize the functional $J_h$ in the sense of \eqref{Exact-Insentitizing}? 
\end{enumerate}}

If $\omega$ and $\Theta$ are strongly included in $\Omega_0$ and $\omega\cap \Theta \not = \emptyset$, {\bf Q1} has been solved in \cite{Lissy-Privat-Simpore}, whereas {\bf Q2} has also been solved in \cite{Lissy-Privat-Simpore} when $\mathcal{E}$ is of dimension $1$ or $2$. The goal of the present article is to extend the results of \cite{Lissy-Privat-Simpore} to more general geometric settings and to the more general questions above-mentioned.
To be more precise, in the next section,  we will distinguish between the cases where $\omega$ and $\Theta$ intersect or not (see Fig.~\ref{fig:scheme} below), since approaches to deal with them and the results obtained are fairly different.  In the case $\omega \cap \Theta = \emptyset$, {\bf Q1} will be tackled in Theorem~\ref{Thm-Approx-Case1} and {\bf Q2} will be tackled in Theorems~\ref{Thm-Case1-Exact-FD} . In the case $\omega \cap \Theta \not = \emptyset$,  {\bf Q1'} will be tackled in Theorem  \ref{Thm-Approx-Case2} and {\bf Q2'} will be tackled in Theorem and \ref{Thm-Case2-Exact-FD}. Finally, we will provide two partial answers to {\bf Q3} in Theorems~\ref{Thm-insensitizing-Boundary} and \ref{Thm-Negative-Theta-Omega}.

\begin{remark}
	{According to Proposition \ref{Prop-J-differentiable}, although all the above questions \emph{a priori} depend on $\xi \in L^2(0,T; L^2(\R^d))$, in reality, they only depend on the restriction on $\xi$ to $\Omega_0$. Thus, in the following, we shall simply take $\xi \in L^2(0,T; L^2(\Omega_0))$ (extended by $0$ on $\R^d$) without loss of generality.}
\end{remark}

\begin{figure}[h!]
\centering
\begin{minipage}{7cm}
\begin{tikzpicture}[scale=1.5]
\coordinate (a) at (0,1);
\coordinate (b) at (3,1);
\coordinate (c) at (3,0);
\coordinate (d) at (1.5,0);
\coordinate (e) at (0,0);
\path[draw,use Hobby shortcut,closed=true]
(a)..(b)..(c)..(d)..(e);

\draw[rotate around={45:(1,0)}]  (1,1) ellipse (0.5 and 0.2);
\draw[rotate around={60:(2,0)}]  (2.5,1) ellipse (0.2 and 0.4);
\draw (0.3,0.3) node[below]{$\Omega_0$} ;
\draw (0.25,0.83) node[below]{$\Theta$} ;
\draw (1.4,1.1) node[below]{$\omega$} ;
\end{tikzpicture}
\end{minipage}
\begin{minipage}{7cm}
\begin{tikzpicture}[scale=1.5]
\coordinate (a) at (0,1);
\coordinate (b) at (3,1);
\coordinate (c) at (3,0);
\coordinate (d) at (1.5,0);
\coordinate (e) at (0,0);
\path[draw,use Hobby shortcut,closed=true]
(a)..(b)..(c)..(d)..(e);

\draw[rotate around={45:(2.2,1)}]  (1.8,1.2) ellipse (0.5 and 0.2);
\draw[rotate around={60:(2,0)}]  (2.5,1) ellipse (0.2 and 0.4);
\draw (0.3,0.3) node[below]{$\Omega_0$} ;
\draw (1.9,1.2) node[below]{$\Theta$} ;
\draw (1.4,1.1) node[below]{$\omega$} ;
\end{tikzpicture}
\end{minipage}
\caption{The two main situations investigated: (left) the intersection set of $\omega$ and $\Theta$ is empty; (right) the intersection set of $\omega$ and $\Theta$ is nonempty. \label{fig:scheme}}
\end{figure}
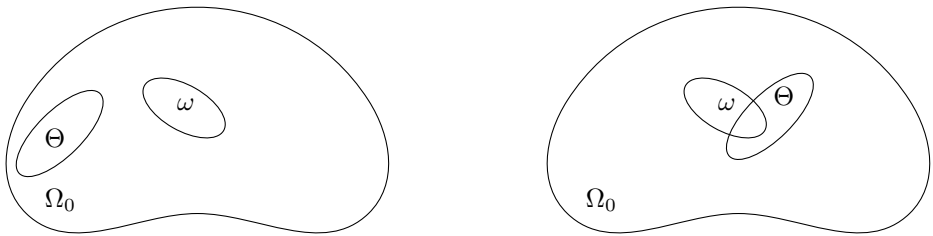

\subsection{Main results}
{As we said,} the \red{aforementioned insensitizing} problems will be strongly dependent of the relative geometry of the various sets $\omega$, $\Theta$ and $\Omega_0$, and in particular depending if the set $\omega \cap \Theta$ is empty or not{, but we will always make the following minimal assumption} on both sets $\omega$ and $\Theta$:
\begin{equation}\tag{$\Hreg$}\label{regMinOmThe}
\text{$\omega$ and $\Theta$ are two nonempty open subsets of $\Omega_0$.}
\end{equation}

\newpage

\begin{center}
\fbox{\textbf{\large Case $\omega \cap \Theta = \emptyset$ and $\Theta \Subset \Omega_0$}}
\end{center} 

To be more precise, the first geometric setting we consider is the following:
\begin{equation}
	\label{Case-1-emptyset}
	\omega\text{ and }\Theta\text{ satisfy }\eqref{regMinOmThe},\quad 
	\omega \cap \Theta = \emptyset, 
	%\quad \Theta \text{ is a \textcolor{blue}{Lipschitz} domain},
	 \quad \Theta \Subset \Omega_0, 
	 \text{ and } \Omega_0 \setminus \overline\Theta \text{ is connected}.
\end{equation}

\ma{In this setting, our first result is dedicated to an approximate controllability property. We will then use it in a crucial way to derive one of the main results of this paper, allowing to partially answer the question Q1.}
%our first result is: 
%
\begin{proposition}
	\label{Prop-Approx-Cont-Case1}
	Assume the geometric setting \eqref{Case-1-emptyset}. Then, given any $(f_1, f_2) \in L^2(0,T; L^2( \partial \Omega_0))^2$, for any $\varepsilon >0$ and $\xi \in L^2(0,T; L^2(\Omega_0))$, there exists a control function $h \in L^2(0,T; L^2(\omega))$ such that the solution $(y,z)$ of 
	\begin{equation}
		\label{Coupled-Heat-Approx}
		\left\{ 
			\begin{array}{rcll}
				\dfrac{\partial y}{\partial t}-\Delta y &=&\xi +h\mathbbm{1}_{\omega
}& \text{ in }(0,T) \times \Omega_0,
				\\ 
				y &=& 0 & \text{ on }(0,T) \times \partial \Omega_0,
				\\ 
				y\left( 0,\cdot\right) &=& 0 & \text{ in }\Omega_0,
				\\
				-\dfrac{\partial z}{\partial t}-\Delta z &=&\mathbbm{1}_{\Theta}y & \text{ in }(0,T) \times \Omega_0,
				\\ 
				z &=& 0 & \text{ on }(0,T) \times \partial \Omega_0,
				\\ 
				z\left( T,\cdot\right) &=& 0 & \text{ in }\Omega_0, 
			\end{array}
		\right.
	\end{equation}
	satisfies 
	\begin{equation}
		\label{Approx-Cont}
		\| \partial_n y - f_1 \|_{L^2(0,T; L^2(\partial\Omega_0))} 
		+ 
		\| \partial_n z - f_2 \|_{L^2(0,T; L^2(\partial\Omega_0))} 
		\leq
		\varepsilon.
	\end{equation}
\end{proposition}

\begin{remark}
\begin{itemize}
\item Looking carefully at the proof of Proposition \ref{Prop-Approx-Cont-Case1}, given in Section \ref{Subsec-Proof-Approx-Case1}, it is easy to figure out that the last condition in \eqref{Case-1-emptyset} can be relaxed {into the following one: $\omega$ intersects every connected component $\mathcal O$ of $\Omega_0 \setminus \overline\Theta$ verifying that $\overline{\mathcal O} \cap \partial \Omega \not = \emptyset$.}

	\item As it is classical for approximate controllability results (see \cite{Zuazua-97,FernandezCaraZuazua1,Ervedoza-Note-2020}), one can reinforce the above results as follows: if $F$ is a finite-dimensional subspace of $(L^2(0, T; L^2(\partial \Omega_0)))^2$ and $\mathbb{P}_F$ denotes the orthogonal projection on $F$ in $(L^2(0, T; L^2(\partial \Omega_0)))^2$, then, for any $(f_1, f_2) \in L^2(0,T; L^2( \partial \Omega_0))^2$, for any $\varepsilon >0$ and $\xi \in L^2(0,T; L^2(\Omega_0))$, there exists a control function $h \in L^2(0,T; L^2(\omega))$ such that the solution $(y,z)$ of \eqref{Coupled-Heat-Approx} satisfies \eqref{Approx-Cont} and 
	\[
		\mathbb{P}_F (\partial_n y, \partial_n z) = \mathbb{P}_F (f_1, f_2).
	\]
	\end{itemize}
\end{remark}

As we will see, the proof of Proposition \ref{Prop-Approx-Cont-Case1} given in Section \ref{Subsec-Proof-Approx-Case1} mainly relies on a unique continuation property for the adjoint operator, which consists of coupled parabolic equations where the coupling coefficients are disjoint from the observation set. This kind of issues is known to be particularly difficult in the case of coupled parabolic systems (see \textit{e.g.} \cite{mousquetaires} for partial results in one space dimension), and  comes naturally when dealing with insensitizing problems. However, to our knowledge, the only works dealing with control and observation sets which do not intersect in this context are \cite{Micu-Ortega-deTeresa} in a $1$d case and \cite{Kavian-De-Teresa}. Though, our result is different, since the unique continuation property we need to prove Proposition \ref{Prop-Approx-Cont-Case1} is not the one in \cite{Micu-Ortega-deTeresa,Kavian-De-Teresa}.

A straightforward application of Proposition \ref{Prop-Approx-Cont-Case1}, proved in Section \ref{Subsec-Proof-ApproxCont-Case1}, is the following one: 
\begin{theorem}
	\label{Thm-Approx-Case1}
	Assume the geometric setting \eqref{Case-1-emptyset}. Then,  for all $\xi \in L^2(0,T; L^2(\Omega_0))$ and $\varepsilon >0$, there exists a control function $h \in L^2(0,T; L^2(\omega))$ such that the solution $(y_0,z_0)$ of \eqref{y-Omega-0-h} satisfies
	\begin{equation}
		\label{Approx-Insensitzing-Equivalent}
		\norm{ \int_0^T \partial_n y_0 \partial_n z_0 \, dt}_{L^1(\partial \Omega_0)} \leq \varepsilon.
	\end{equation}
	In other words, according to \eqref{Diff-J-V-2}, the functional $J_h$ is $\varepsilon$-approximately insensitized by $h \in L^2(0,T; L^2(\omega))$ at $\Omega_0$ in the sense of \eqref{Approx-Insensitizing}.
\end{theorem}
%
%\begin{remark}
%Because of the expression of $\left.\dfrac{d }{d \tau}\left({J}_h\left(\Omega_{\tau\mathbf{V}} \right)\right)\right|_{\tau=0}$ given by \eqref{Diff-J-V-2}, it is clear that \eqref{Approx-Insensitizing} is a consequence of \eqref{Approx-Insensitzing-Equivalent}, by using the fact that $\Vert \mathbf{V}\Vert_{L^\infty}\leq \Vert \mathbf{V}\Vert_{W^{3,\infty}}$.
%\end{remark}

One can actually prove that the functional $J_h$ can be made exactly insensitized to any finite-dimensional vector space of $W^{3, \infty}(\R^d, \R^d)$:
\begin{theorem}
	\label{Thm-Case1-Exact-FD}
	Assume the geometric setting \eqref{Case-1-emptyset}. Let  $\mathcal{E}$ be a finite-dimensional linear subspace of $W^{3, \infty}(\R^d, \R^d)$. Then,  for all $\xi \in L^2(0,T; L^2(\Omega_0))$ and for all $\varepsilon>0$, there exists a control $h\in L^2(0,T; L^2(\omega))$ that insensitizes $J_h$ exactly for $\mathcal{E}$ in the sense of \eqref{Finite-Insensitizing} and that $\varepsilon$-approximately insensitizes $J_h$ in the sense of \eqref{Approx-Insensitizing}.
\end{theorem}
%

%These two results are the first ones obtained in the context of a control set $\omega$ and an observation set $\Theta$ which do not intersect. 

In fact, the main difficulty in Theorem \ref{Thm-Case1-Exact-FD} is the construction of $h$ such that $J_h$ is exactly insensitized for $\mathcal{E}$ in the sense of \eqref{Finite-Insensitizing}, since the map 
\begin{equation}
	\label{Derivative-h}
	h \in L^2(0,T; L^2(\omega)) \mapsto  \left(  \int_0^T \partial_n y_0 \partial_n z_0 \, dt\right) \in  L^1(\partial \Omega_0), 
\end{equation}
where the pair $(y_0,z_0)$ solves \eqref{y-Omega-0-h}, is not linear in $h$ even in the case $\xi = 0$, but \red{quadratic}. Therefore, we use techniques specifically designed to deal with such kind of non-linearities, which consists in choosing the control functions in a vector space of much larger dimension than the number of constraints. Similarly to what has been done in another context for the stabilizability of the Navier-Stokes equation, see \cite{Chowdhury-Erv-2019}, if there are $N$ constraints imposed by the exact insensitizing for $\mathcal E$, we look for control functions in a vector space of size (at most) $2N$ which is suitably designed. In particular, even if there is one constraint (\textit{i.e.} if $\mathcal{E}$ is a vector space of dimension $1$), we look for the control function in a vector space of dimension (at most) $2$, thus preventing possible obstructions that may appear due to the quadratic nature of the map in \eqref{Derivative-h} (see \textit{e.g.} \cite{Beauchard-Marbach-2017}). 
Details of the proof are given in Section~\ref{Subsec-Proof-Exact-Case1}.

\begin{center}
\fbox{\textbf{\large Case $\omega \cap \Theta \neq \emptyset$}}
\end{center}  
The second geometric setting we consider is the case 
\begin{equation}
	\label{Case-2-non-emptyset}
	\omega\text{ and }\Theta\text{ satisfy }\eqref{regMinOmThe},\quad 	\omega \cap \Theta \neq \emptyset. 
\end{equation}
This case is more favorable since the control set $\omega$ meets the observation set \red{$\Theta$}, and as we shall see afterwards, not only do all previously established results remain true, but this also allows to prove the existence of even better insensitizing controls. 

\paragraph{On the $\varepsilon$-approximate insensitizing problem.}

 To start with, we first claim that Proposition \ref{Prop-Approx-Cont-Case1} can be reinforced under the geometric setting \eqref{Case-2-non-emptyset}. 

\begin{proposition}
	\label{Prop-Approx-Cont-Case2}
	Assume the geometric setting \eqref{Case-2-non-emptyset}. Then, given any $(f_1, f_2) \in (L^2(0,T; L^2( \partial \Omega_0)))^2$,  $y_T \in L^2(\Omega)$, any $\varepsilon >0$ and $\xi \in L^2(0,T; L^2(\Omega_0))$, there exists a control function $h \in L^2(0,T; L^2(\omega))$ such that the solution $(y,z)$ of \eqref{Coupled-Heat-Approx} satisfies 
	\begin{equation}
		\label{Approx-Cont-Case2}
		\| \partial_n y - f_1 \|_{L^2(0,T; L^2(\partial\Omega_0))} 
		+ 
		\| \partial_n z - f_2 \|_{L^2(0,T; L^2(\partial\Omega_0))} 
		+
		\| y(T) - y_T \|_{L^2(\Omega_0)}
		\leq
		\varepsilon.
	\end{equation}

	Besides, if {the source term} $\xi \in L^2(0,T; L^2(\Omega_0))$ is null-controllable in the sense that there exists $h_{nc} \in L^2(0,T; L^2(\omega))$ such that the solution $y_{nc}$ of
	\begin{equation}
		\label{y-nc}
		\left\{ 
			\begin{array}{rcll}
				\dfrac{\partial y_{nc}}{\partial t}-\Delta y_{nc} &=&\xi +h_{nc} \mathbbm{1}_{\omega
}& \text{ in }(0,T) \times \Omega_0,
				\\ 
				y_{nc} &=& 0 & \text{ on }(0,T) \times \partial \Omega_0,
				\\ 
				y_{nc}\left( 0,\cdot\right) &=& 0 & \text{ in }\Omega_0,
			\end{array}
		\right.
	\end{equation}
	 satisfies 
	\begin{equation}
		\label{NC-Assumption}
		y_{nc}(T) = 0 \hbox{ in } \Omega_0, 
	\end{equation}
	then, given any $(f_1, f_2) \in L^2(0,T; L^2( \partial \Omega_0))^2$, for any $\varepsilon >0$, there exists a control function $h \in L^2(0,T; L^2(\omega))$ such that the solution $(y,z)$ of \eqref{Coupled-Heat-Approx} satisfies
		\begin{equation}
		\label{Approx-Cont-2}
		\| \partial_n y - f_1 \|_{L^2(0,T; L^2(\partial\Omega_0))} 
		+ 
		\| \partial_n z - f_2 \|_{L^2(0,T; L^2(\partial\Omega_0))} 
		\leq
		\varepsilon.
	\end{equation}
	and 
	\begin{equation}
		\label{NullControl}
		y(T) = 0 \text{ in } \Omega_0.
	\end{equation}
\end{proposition}
\begin{remark}
Determining a control function $h$ such that the solution $y_{nc}$ of \eqref{y-nc} satisfies \eqref{NC-Assumption} is the well-known null-controllability problem with source term for the heat equation. This issue has been much investigated. By using duality arguments, this issue can be recast in terms of a so-called ``observability inequality''. More precisely, one can deduce from \cite[Lemma 2.1]{FursikovImanuvilov} that there exists $C>0$ such that 
\begin{equation}\label{fc}
	\left\| \exp\left( - \frac{C}{T-t}\right)\varphi(t) \right\|_ {L^2(0,T; L^2(\Omega_0))} \leq C \| \varphi \|_{L^2(0,T; L^2(\omega))},
\end{equation}
 where $\varphi$ denotes the solution of the backward adjoint system
	\begin{equation*}
		\left\{ 
			\begin{array}{rcll}
				-\dfrac{\partial \varphi}{\partial t}-\Delta \varphi &=&0 & \text{ in }(0,T) \times \Omega_0,
				\\ 
				\varphi &=& 0 & \text{ on }(0,T) \times \partial \Omega_0,
				\\ 
				\varphi\left( T,\cdot\right) &=& \varphi^T & \text{ in }\Omega_0,
			\end{array}
		\right.
	\end{equation*}
From the observability inequality \eqref{fc}, we can deduce that  for any $\xi \in L^2(0,T; L^2(\Omega_0))$ with $e^{\frac{C}{T-t}}\xi \in L^2(0,T; L^2(\Omega_0))$, one can find a null control $h_{nc}$ to \eqref{y-nc}, \textit{i.e.} a function $h_{nc} \in L^2(0,T; L^2(\omega))$ such that the solution $y_{nc}$ of \eqref{y-nc} satisfies \eqref{NC-Assumption}, as explained in \cite[Proof of Theorem 2.1]{FursikovImanuvilov}.

\end{remark}

Again, Proposition \ref{Prop-Approx-Cont-Case2} is based on suitable unique continuation properties for the adjoint equation. However, here, since $\omega \cap \Theta \neq\emptyset$, the arguments are more standard for the proof of \eqref{Approx-Cont-Case2} {than for the proof of Proposition \ref{Prop-Approx-Cont-Case1}}. The possibility of further imposing \eqref{NullControl} when $\xi$ is a source term that can be null-controlled is much more subtle, and amounts to a suitable use of duality arguments, inspired by \cite{Lions-1992}, and of observability estimates for the heat equation given in \cite{FursikovImanuvilov}. Details of the proof are given in Section~\ref{Subsec-Proof-ApproxCont-Case2}.

As before, a straightforward application of Proposition \ref{Prop-Approx-Cont-Case2} is the following result, whose proof is postponed to Section~\ref{Subsec-Proof-Approx-Case2}.
\begin{theorem}
	\label{Thm-Approx-Case2}
	Assume the geometric setting \eqref{Case-2-non-emptyset}. Then, for all $\xi \in L^2(0,T; L^2(\Omega_0))$, $y_T \in L^2(\Omega_0)$ and $\varepsilon >0$, there exists a control function $h \in L^2(0,T; L^2(\omega))$ such that the solution $(y_0,z_0)$ of \eqref{y-Omega-0-h} satisfies \eqref{Approx-Insensitzing-Equivalent} and 
	\begin{equation}
		\label{Approx-Cont-T}
		\norm{y_0(T) - y_T}_{L^2(\Omega_0)} \leq \varepsilon.
	\end{equation}
	In other words, according to \eqref{Diff-J-V-2}, the functional $J_h$ is $\varepsilon$-approximately insensitized by a control $h \in L^2(0,T; L^2(\omega))$ at $\Omega_0$ in the sense of \eqref{Approx-Insensitizing}. \red{Furthermore, the function $h$ also $\varepsilon$-approximately controls the state $y_0$ of \eqref{y-Omega-0-h} at time $T$, in the sense that \eqref{Approx-Cont-T} is verified.}

	Similarly,  if {the source term} $\xi \in L^2(0,T; L^2(\Omega_0))$ is null-controllable in the sense that there exists $h_{nc} \in L^2(0,T; L^2(\omega))$ such that the solution $y_{nc}$ of \eqref{y-nc} satisfies \eqref{NC-Assumption}, then, there exists a control function $h$ such that solution $(y_0,z_0)$ of \eqref{y-Omega-0-h} satisfies \eqref{Approx-Insensitzing-Equivalent} and 
	\begin{equation}
		\label{NC-and-Approx-Insenstizing}
		y_0(T) = 0 \text{ in } \Omega_0.
	\end{equation}
	In other words, if {the source term} $\xi$ is null-controllable at time $T>0$, the functional $J_h$ is $\varepsilon$-approximately insensitized by a control $h \in L^2(0,T; L^2(\omega))$ at $\Omega_0$ in the sense of \eqref{Approx-Insensitizing}, which also steers the state $y_0$ of \eqref{y-Omega-0-h} exactly to $0$ at time $T$. 
\end{theorem}

One can also improve Theorem \ref{Thm-Case1-Exact-FD} in the case of the geometric setting \eqref{Case-2-non-emptyset}:
\begin{theorem}
	\label{Thm-Case2-Exact-FD}
	Assume the geometric setting \eqref{Case-2-non-emptyset}, and let $\mathcal{E}$ be a finite-dimensional subspace of $W^{3, \infty}(\R^d, \R^d)$. 
	
	Then, for all $ \xi \in L^2(0,T; L^2(\Omega_0))$ and $y_T \in L^2(\Omega_0)$, for all $\varepsilon>0$ , there exists a control $h\in L^2(0,T; L^2(\omega))$ that insensitizes $J_h$ exactly for $\mathcal{E}$ in the sense of \eqref{Finite-Insensitizing}, $\varepsilon$-approximately insensitizes $J_h$ in the sense of \eqref{Approx-Insensitizing}, and which approximately controls $y_0$ at time $T$ in the sense of \eqref{Approx-Cont-T}.
	
	Besides, if {the source term} $ \xi \in L^2(0,T; L^2(\Omega_0))$ is null-controllable, then, there exists a control $h\in L^2(0,T; L^2(\omega))$ that insensitizes $J_h$ exactly for $\mathcal{E}$ in the sense of \eqref{Finite-Insensitizing}, $\varepsilon$-approximately insensitizes $J_h$ in the sense of \eqref{Approx-Insensitizing}, and which steers $y_0$ to $0$ at time $T$ in the sense of \eqref{NC-and-Approx-Insenstizing}.

\end{theorem}
Here again, the proof of Theorem \ref{Thm-Case2-Exact-FD} is a rather simple adaptation of the one of Theorem \ref{Thm-Case1-Exact-FD}, based on the stronger results given by Proposition \ref{Prop-Approx-Cont-Case2}.
\begin{remark}
Remark that Theorems \ref{Thm-Approx-Case2}  and \ref{Thm-Case2-Exact-FD} can be reinterpreted in terms of \textit{robustness} \red{of null and approximate controllability properties}: they notably enable us to find a null or approximate control $h$ for the heat equation 
	\begin{equation*}
		\left\{ 
			\begin{array}{rcll}
				\dfrac{\partial y}{\partial t}-\Delta y &=&\xi +h\mathbbm{1}_{\omega
}& \text{ in }(0,T) \times \Omega_0,
				\\ 
				y &=& 0 & \text{ on }(0,T) \times \partial \Omega_0,
				\\ 
				y\left( 0,\cdot\right) &=& 0 & \text{ in }\Omega_0,
			\end{array}
		\right.
	\end{equation*}
	so that the functional $J_h$ is robust to small variations of the boundary, in the sense that this control makes $J_h$ insensitive at the first order to small perturbations of the boundary. 
\end{remark}

\begin{center}
\fbox{\textbf{\large On the exact insensitizing problem.}}
\end{center}  

%\paragraph{On the exact insensitizing problem.}
%
Note that in both geometric settings discussed so far, the question of exact insensitizing has not been addressed. We now propose to study some cases in which we can solve the insensitizing problem. Let us start with the rather straightforward case $ \Theta \Subset \omega$. 
\begin{proposition}
	\label{Prop-Cas-Facile}
	Let $\omega$ and $\Theta$ be non-empty open subsets of $\Omega_0$ such that
	\begin{equation}
		\label{Assumption-EasyCase}
		\Theta \Subset \omega.
	\end{equation}
	Then, for all $\xi \in L^2(0,T; L^2(\Omega_0))$, there exists $h \in L^2(0,T; L^2(\omega))$ such that the functional $J_h$ in \eqref{def-J} is exactly insensitized in the sense of \eqref{Exact-Insentitizing}.
\end{proposition}
Proposition \ref{Prop-Cas-Facile}, proved in Section \ref{Sec-Proof-Prop-Cas-Facile}, in fact considers an easy case, in which we can ensure that with a suitable choice of a control function, the solution $y_0$ of \eqref{y-Omega-0-h} vanishes in $(0,T) \times \Theta$, so that  the associated function $z_0$ {satisfying}  \eqref{y-Omega-0-h} vanishes in $(0,T) \times \Omega_0$ and the result easily follows from \eqref{Diff-J-V-2}.

\medskip

Let us now consider a more subtle case, in which the outer boundary of $\Theta$ is included in $\omega$.
\begin{theorem}
	\label{Thm-insensitizing-Boundary}
	Let $\omega$ and $\Theta$ be smooth non-empty open subsets of $\Omega_0$ such that
	\begin{equation}
		\label{Assumption-insensitizing-Boundary}
		\partial \Theta \text{ has only one connected component, } \quad \Theta \Subset  \Omega_0, \quad \text{ and } \quad \partial \Theta \subset \omega.
	\end{equation}
	Then, for all $\xi \in L^2(0,T; L^2(\Omega_0))$, there exists $h \in L^2(0,T; L^2(\omega))$ such that the functional $J_h$ in \eqref{def-J} is exactly insensitized in the sense of \eqref{Exact-Insentitizing}.
\end{theorem}
The strategy to prove this theorem is to choose the control function $h$ such that the solution $z_0$ of \eqref{y-Omega-0-h} vanishes close to the boundary $\partial \Omega_0$. Thus, using \eqref{Diff-J-V-2}, the exact insensitizing of $J_h$ will immediately follow. In order to do that, we will interpret the function $y_0$ in \eqref{y-Omega-0-h} as a control function for $z_0$ whose goal is to impose the condition $z_0 = 0$ outside $(0,T) \times \Theta$, and we then define $h$ in terms of $y$ by \eqref{y-Omega-0-h}. We refer to Section \ref{Sec-Proof-Thm-insensitizing-Boundary} for the proof of Theorem \ref{Thm-insensitizing-Boundary}.
%

%One can easily check by inspecting the proof of Theorem \ref{Thm-insensitizing-Boundary} that the assumption \eqref{Assumption-insensitizing-Boundary} can be easily weakened to the following one: $\Theta \Subset  \Omega$ and there exist a connected component $\Gamma_\Theta$ of $\partial \Theta$ and an open subset $\tilde\Theta$ with boundary $\partial\tilde\Theta = \Gamma_\Theta$ such that $\tilde \Theta \Subset  \Omega$ and $\Gamma_{\Theta} \subset \omega$. This is in fact a fancy way to express that the outer boundary of $\Theta$ is contained in $\omega$. 

These two positive results should very likely not be considered as the usual case. In fact, we can discuss the case $\Theta = \Omega_0$ with more details:
\begin{theorem}
	\label{Thm-Negative-Theta-Omega}
	Assume that {$\Omega_0$  is smooth (of class $\mathscr{C}^\infty$)}, that $\Theta = \Omega_0$ and that $\omega$ is a non-empty open subset of $\Omega_0$ such that $\omega \Subset \Omega_0$. Then, there exists a function $\xi \in L^2(0,T; L^2(\Omega_0))$ such that there exists no $h \in L^2(0,T; L^2(\omega))$ such that $J_h$ satisfies \eqref{Exact-Insentitizing}. In other words, there are some $\xi \in L^2(0,T; L^2(\Omega_0))$ such that the exact insensitizing problem cannot be solved. 
\end{theorem}
The proof of this result is given in Section \ref{Sec-Proof-Thm-Negative-Theta-Omega}, and in fact only involves regularity issues.

\subsection{Bibliographical comments}\label{sec:SoA}

We comment briefly on the bibliography, emphasizing particularly the works related to the heat equation (linear or non-linear) and dedicated approaches to solving problems related to {functional insensitizing}.

The question of {functional insensitizing} has been first introduced in \cite{Lions-Sentinelles-1992,Bodart-Fabre}. However, in \cite{Lions-Sentinelles-1992,Bodart-Fabre}, the functionals under consideration were insensitized with respect to perturbation of the initial datum or of the source term, while we are discussing a new kind of insensitizing, with respect to perturbation of the boundary. 

Still, our approach is of course strongly inspired by the one developed in \cite{Lions-Sentinelles-1992,Bodart-Fabre}, in which it was shown how unique continuation properties can be used to solve approximate insensitizing problems. It was then further developed to many settings, in particular when the control set and the observation set intersect. 

Regarding the standard issue of insensitizing of a given functional (often the $L^2$ norm of the state in some observation subset) involving the solution of the heat equation with respect to initial data, the general approach consists in recasting the (exact or approximated) insensitizing property in terms of an adjoint state, leading to consider a coupled system of forward-backward heat equations. Hence, exact insensitizing comes to investigate a null-controllability property which can in general be recasted through an observability inequality (see \cite{deTeresa-00,Bodart-GB,Boyer-HS-dT} where Carleman based approaches are considered and \cite{dT-Zuazua-09}, in which a Fourier approach is used).  
The question of $\varepsilon$-approximated insensitizing comes in general to solve an approximate controllability problem, leading to derive a unique continuation property (see \cite{Micu-Ortega-deTeresa,Kavian-De-Teresa}, in which spectral methods are employed).
We also mention \cite{Guerrero,Gueye,Carreno-Gueye,Carreno-Guer-Gueye,Carreno} where a functional involving the solution of another equations arising in Fluid Mechanics is considered. 

\subsection{{Further comments and open problems}}

In this article, we investigate and discuss three insensitizing properties with respect to boundary variations. To conclude this introduction, we outline {three} open issues and hints that complement the study presented in this article and that we plan to address in the future.

\paragraph{Open problem \#1.}
{Note that we were not able to answer questions  {\bf Q1'} and  {\bf Q2'} when $\omega \cap \Theta = \emptyset$. The main difference with the case $\omega \cap \Theta \neq \emptyset$ is that the approximate controllability results we are able to prove in the case $\omega \cap \Theta = \emptyset$ is weaker than in the case $\omega \cap \Theta \neq \emptyset$, compare Proposition \ref{Prop-Approx-Cont-Case1} and Proposition \ref{Prop-Approx-Cont-Case2}. As one can check from the proofs, the stronger statement in Proposition \ref{Prop-Approx-Cont-Case2} comes by duality from unique continuation properties for a coupled parabolic system, namely the unique continuation property \eqref{UC-Adj-Coupled-Case2} for the solutions of \eqref{Coupled-Heat-Approx-Adj-Case2}. Whether this unique continuation property holds when $\omega\cap\Theta =\emptyset$ is an open problem.
}

\paragraph{Open problem \#{2}.} Can we answer the three questions {{\bf Q1} -- {\bf Q2} -- {\bf Q3}} posed in this article when the statement of problem is modified as follows:
\begin{itemize}
\item the heat equation \eqref{y-Omega-h} is replaced by more general controlled parabolic equations, e.g. semi-linear problems, Stokes or Navier-Stokes systems \ma{(the case of insensitizing with respect to initial data in this last case has been notably investigated in \cite{Guerrero})} or a controlled wave equation \ma{(the case of insensitizing with respect to initial data in this last case has been notably investigated in \cite{Tebou1,Tebou2,Alabau-Boussouira})}. Note that this is very likely difficult to handle with the arguments developed here, since they are based on approximate controllability statements (recall Propositions \ref{Prop-Approx-Cont-Case1} and \ref{Prop-Approx-Cont-Case2});

\item the shape functional with respect to which insensitizing is performed is replaced by a more general one of the kind 
$$
\int_0^T\int_{\Theta}j(y(t,x),\nabla y(t,x))\, dxdt,
$$  
where $j:\R^{1+d} \to \R$ is a given function, and $y$ denotes the solution of the considered controlled system.
\end{itemize}

\paragraph{Open problem \#{3}.} Can the answers provided in this article related to exact insensitizing be completed? In particular, what can be expected in the case where $\Theta \subsetneq \Omega_0$? Is it possible to answer positively or negatively to questions {\bf Q1'} and {\bf Q2'} when $\omega \cap \Theta= \emptyset$? Can one identify the set of functions $\xi$ in $L^2(0,T; L^2(\R^d))$ for which {\bf Q3} holds true?

\subsection{Outline}

This article is organized as follows. Section \ref{Sec-Omega-Theta-Empty} studies the case \eqref{Case-1-emptyset}, \textit{i.e.} $\omega \cap \Theta = \emptyset$ and $\Theta \Subset \Omega_0$, and gives the proofs of Proposition \ref{Prop-Approx-Cont-Case1}, Theorem \ref{Thm-Approx-Case1}, and Theorem \ref{Thm-Case1-Exact-FD}. Section \ref{Sec-Omega-Theta-Not-Empty} then focuses on the case \eqref{Case-2-non-emptyset}, \textit{i.e.} $\omega \cap \Theta \neq \emptyset$, and provides the proofs of Proposition \ref{Prop-Approx-Cont-Case2}, Theorem \ref{Thm-Approx-Case2} and Theorem \ref{Thm-Case2-Exact-FD}. Section \ref{Sec-Exact-insensitizing} then presents the proofs of the results on exact insensitizing, namely Proposition \ref{Prop-Cas-Facile}, Theorem \ref{Thm-insensitizing-Boundary} and Theorem \ref{Thm-Negative-Theta-Omega}.
%
%
%%%%%%%%%%%%%%%%%%%%%%%%%%%%%%%%%%%%%%%%%%%%%%%%%%%%%%%%%%%%%%%%%%%%%%%%%%%%%%%%%%%%%%%%%%%%%%%%%%%%%%%%%%%%%%%%%%%%%%%%%%%%%%%%%%%%%%%
%

\section{The case $\omega \cap \Theta = \emptyset$ and $\Theta \Subset \Omega_0$.}
\label{Sec-Omega-Theta-Empty} 
In this whole section, we assume the geometric setting described in \eqref{Case-1-emptyset}.%, i.e. $\omega \cap \Theta= \emptyset$ and $\Theta \Subset \Omega_0$.

\subsection{Proof of Proposition \ref{Prop-Approx-Cont-Case1}}\label{Subsec-Proof-ApproxCont-Case1}
Proposition \ref{Prop-Approx-Cont-Case1} can be recast in an abstract form into the problem: show that 
\[
	\overline{ \text{Ran\,} \mathcal{L} } = (L^2(0,T; L^2(\partial\Omega_0)))^2, 
\]
where $\mathcal{L}$ is the operator defined for $h \in L^2(0,T; L^2(\omega))$ by $\mathcal{L} h = (\partial_n y_h, \partial_n z_h)\in L^2((0,T)\times\partial\Omega_0)^2$, where $(y_h,z_h)$ solves 	
\begin{equation}\label{csanss}
		\left\{ 
			\begin{array}{rcll}
				\dfrac{\partial y_h}{\partial t}-\Delta y_h &=&h\mathbbm{1}_{\omega
}& \text{ in }(0,T) \times \Omega_0,
				\\ 
				y_h &=& 0 & \text{ on }(0,T) \times \partial \Omega_0,
				\\ 
				y_h\left( 0,\cdot\right) &=& 0 & \text{ in }\Omega_0,
				\\
				-\dfrac{\partial z_h}{\partial t}-\Delta z_h &=&\mathbbm{1}_{\Theta}y_h & \text{ in }(0,T) \times \Omega_0,
				\\ 
				z_h &=& 0 & \text{ on }(0,T) \times \partial \Omega_0,
				\\ 
				z_h\left( T,\cdot\right) &=& 0 & \text{ in }\Omega_0.
			\end{array}
		\right.
	\end{equation}
Notice that $\mathcal L$ is bounded, thanks to \cite[Theorem 5, Page 382]{Evans} and the continuity of the operator $f\in H^2(\Omega_0) \mapsto \partial_n f\in L^2((0,T)\times\partial\Omega_0)$.
Therefore, by standard arguments from functional analysis, Proposition \ref{Prop-Approx-Cont-Case1} is equivalent to showing \red{the dual property} $\text{Ker}\,\mathcal{L}^* = \{0\}$ \red{(we refer for instance to \cite[Corollary~2.18]{Brezis}). Using the arguments developed in \cite[Proposition 2.1]{Lissy-Privat-Simpore}, one can compute explicitly $\mathcal L^*$ and deduce that $\text{Ker}\,\mathcal{L}^* = \{0\}$} is equivalent to the following unique continuation problem: If $( g_1, g_2) \in (L^2(0,T; L^2(\partial\Omega_0)))^2$, and $( \psi, \varphi)$ solves 
\begin{equation}
	\label{Coupled-Heat-Approx-Adj}
		\left\{ 
			\begin{array}{rcll}
				-\dfrac{\partial \psi}{\partial t}-\Delta \psi &=&\mathbbm{1}_{\Theta}\varphi & \text{ in }(0,T) \times \Omega_0,
				\\ 
				\psi &=& g_1 & \text{ on }(0,T) \times \partial \Omega_0,
				\\ 
				\psi \left( T,\cdot\right) &=& 0 & \text{ in }\Omega_0, 
				\\
				 \dfrac{\partial \varphi}{\partial t}-\Delta \varphi &=& 0 & \text{ in }(0,T) \times \Omega_0,
				\\ 
				\varphi &=& g_2 & \text{ on }(0,T) \times \partial \Omega_0,
				\\ 
				\varphi\left( 0,\cdot\right) &=& 0 & \text{ in }\Omega_0,
			\end{array}
		\right.
\end{equation}
then, we have the following unique continuation property:
\begin{equation}
	\label{UC-Adj-Coupled}
	\psi = 0 \text{ in } (0,T) \times \omega \Rightarrow g_1 =  g_2 = 0.
\end{equation} 

We now prove this unique continuation property. Let $(g_1, g_2) \in (L^2(0,T; L^2(\partial\Omega_0)))^2$ be such that the solution $(\psi, \varphi)$ of \eqref{Coupled-Heat-Approx-Adj} satisfies $\psi = 0$ in $(0,T) \times \omega$. From \cite[(15.17-18), Page 86]{LionsMagenes2}, we notably infer that $(\psi,\varphi)\in L^2((0,T)\times\Omega_0)$.

We first work in the set $(0,T) \times (\Omega_0 \setminus \overline{\Theta})$. There, $\psi$ satisfies the usual backward heat equation 
$$-\dfrac{\partial \psi}{\partial t}-\Delta \psi =0\text{ in }(0,T) \times(\Omega_0 \setminus \overline{\Theta}),$$
where we do not specify any initial or boundary conditions.
 Therefore, since $\omega \subset \Omega_0\setminus {\overline{\Theta}}$ and $\Omega_0\setminus {\overline{\Theta}}$ is connected thanks to Assumption \eqref{Case-1-emptyset}, \red{
using the Holmgren theorem for the heat equation (see \textit{e.g.} \cite[Section~8.6]{Hormander-2003}), we thus infer that $\psi = 0$ in $(0,T) \times (\Omega_0 \setminus \overline{\Theta})$. 
 In particular, using \cite[Proposition 7.1.3]{TWBook}, since $\psi$ the unique weak solution of \eqref{Coupled-Heat-Approx-Adj} in the sense of \cite[Definition 7.1.2, Page 342]{TWBook},  we have that for any $\zeta\in H^2(\Omega_0)\cap H^1_0(\Omega_0)$ such that $\Delta \zeta \in H^1_0(\Omega_0)$, for any $t\in [0,T]$, we have 
 $$\int_0^t\int_{(0,T)\times\partial\Omega_0}g_1\partial_n \zeta d\sigma dt =0.$$
 Differentiating this inequality with respect to $t$ (which is possible because $g_1(\cdot,x)\partial_n \zeta(x) \in L^2(0,T)$ for almost all $x\in (0,T)\times\partial\Omega_0$, since the Neumann trace operator  $f\in H^2(\Omega_0) \mapsto \partial_n f\in L^2(\partial\Omega_0)$ is well-defined), we obtain that for any $t\in [0,T]$, we have 
$$ \int_{\partial\Omega_0}g_1(t,\cdot)\partial_n \zeta d\sigma =0.$$
Since $\{\zeta \in H^2(\Omega_0)\cap H^1_0(\Omega_0)| \Delta \zeta \in H^1_0(\Omega_0)\}$ is dense in $H^2(\Omega_0)\cap H^1_0(\Omega_0)$ for the $H^2$-norm, using the continuity of the Neumann trace operator  $f\in H^2(\Omega_0) \mapsto \partial_n f\in L^2((0,T)\times\partial\Omega_0)$), we deduce that for any $\zeta\in H^2(\Omega_0)\cap H^1_0(\Omega_0)$ , we have, for any $t\in [0,T]$,
$$ \int_{\partial\Omega_0}g_1(t,\cdot)\partial_n \zeta d\sigma =0.$$
Since the extended trace operator $f\in H^2(\Omega)\mapsto (f_{|\partial\Omega_0},\partial_n f_{|\partial\Omega_0}) \in H^{3/2}(\partial\Omega_0)\times  H^{1/2}(\partial\Omega_0) $ is surjective, taking into account that $0\in H^{3/2}(\partial\Omega_0)$, we deduce that for any $h\in H^{1/2}(\partial\Omega_0)$, we have, for any $t\in [0,T]$,
$$ \int_{\partial\Omega_0}g_1(t,\cdot) h d\sigma =0.$$
Since  $H^{1/2}(\partial\Omega_0)$ is dense in $L^2(\partial\Omega_0)$ for the $L^2(\partial\Omega_0)$-norm, we deduce that  $g_1 = 0$ on $(0,T)\times \partial\Omega_0$,  so that $\psi \equiv 0$ on $\Omega_0\setminus \overline{\Theta}$ in the classical sense. Notably, $\partial_n \psi= 0$ on the whole boundary $\partial \Omega_0$ since $\Theta \Subset \Omega_0$.}
  
 \red{We next remark that, by standard local regularity results for solutions of the heat equation, the function $\varphi$ is smooth away from the boundary $\{0,T\}\times\partial \Omega_0$ (of class $\mathscr{C}^\infty((0,T) \times \Omega_0)$). Accordingly, $\mathbbm{1}_{\Theta}\varphi $ belongs to $L^2(0,T; L^2(\Omega))$, and by standard maximal regularity results, the solution $\psi$ belongs to $L^2(0,T; H^2_{loc}(\Omega_0)) \cap H^1(0,T; L^2_{loc}(\Omega_0))$.} 
%
%We next consider the equation \eqref{Coupled-Heat-Approx} in $(0,T) \times \Theta$. There, we have 
%%
%\begin{equation}
%	\label{Coupled-Heat-Approx-Adj-Theta}
%		\left\{ 
%			\begin{array}{rcll}
%				-\dfrac{\partial \psi}{\partial t}-\Delta \psi &=& \varphi & \text{ in }(0,T) \times \Theta,
%				\\ 
%				\psi &=& 0 & \text{ on }(0,T) \times \partial \Theta,
%				\\ 
%				\partial_n \psi &=& 0 & \text{ on }(0,T) \times \partial \Theta,
%				\\ 
%				\psi \left( T,\cdot\right) &=& 0 & \text{ in }\Theta, 
%				\\
%				 \dfrac{\partial \varphi}{\partial t}-\Delta \varphi &=& 0 & \text{ in }(0,T) \times \Theta,
%				\\
%				\varphi\left( 0,\cdot\right) &=& 0 & \text{ in }\Theta,
%			\end{array}
%		\right.
%\end{equation}
%
\red{In view of this regularity, and due to the fact that $\psi = 0$ in a neighbourhood of $(0,T) \times \partial \Omega_0$, we can multiply equation \eqref{Coupled-Heat-Approx-Adj}$_{(1)}$  by $\varphi$ and perform several integration by parts, using $g_1 = \partial_n\psi =0$ on $\partial\Omega_0$, and $\psi = 0$ in $(0,T) \times (\Omega_0 \setminus \overline{\Theta})$:} 
\begin{align*}
	& \int_0^T \int_{\Theta} |\varphi|^2 \, dx dt 
	= 
	\int_0^T \int_{\Omega_0} \varphi \left( - \dfrac{\partial \psi}{\partial t}-\Delta \psi\right) \, dx dt
	\\
	&
	= 
	- \left.\int_{\Omega_0} \varphi(\cdot,x) \psi(\cdot, x)\, dx\right|_0^T
	- 
	\int_0^T \int_{\partial\Omega_0} \varphi \partial_n \psi d\sigma dt 
	+
	\int_0^T \int_{\partial\Omega_0} \partial_n \varphi  \psi d\sigma dt 
	+
	\int_0^T \int_{\Omega_0}  \left(  \dfrac{\partial \varphi}{\partial t}-\Delta \varphi \right) \psi \, dx dt
%	\\
%	& = \int_0^T \int_{\Theta}  \left(  \dfrac{\partial \varphi}{\partial t}-\Delta \varphi \right) \psi \, dx dt
	\\&=0.
\end{align*}
%\textcolor{blue}{The last line comes from the fact that $\psi = 0$ in $(0,T) \times (\Omega_0 \setminus \overline{\Theta})$.}
Therefore, $\varphi = 0$ in $(0,T) \times \Theta$, and by the classical unique continuation properties for the heat equation, $\varphi = 0$ in $(0,T) \times \Omega_0$, and in particular $g_2 = 0$. This concludes the proof of \eqref{UC-Adj-Coupled} for the solutions of \eqref{Coupled-Heat-Approx-Adj}, \red{which proves that $\text{Ker}\,\mathcal{L}^* = \{0\}$. Hence, }Proposition \ref{Prop-Approx-Cont-Case1} follows.

\subsection{Proof of Theorem \ref{Thm-Approx-Case1}}\label{Subsec-Proof-Approx-Case1}

The proof of Theorem \ref{Thm-Approx-Case1} mainly reduces to Proposition \ref{Prop-Approx-Cont-Case1}. Indeed, from Proposition \ref{Prop-Approx-Cont-Case1} with $f_1=f_2=0$, for any $\varepsilon >0$, there exists a control function $h \in L^2(0,T; L^2(\omega))$ such that the solution $(y_0,z_0)$ of \eqref{y-Omega-0-h} satisfies
\begin{equation}
	\label{Small-Neumann-Data}
	\norm{\partial_n y_0}_{L^2(0,T;L^2(\partial\Omega_0))} + \norm{\partial_n z_0}_{L^2(0,T;L^2(\partial\Omega_0))}\leq \sqrt{\varepsilon}. 
\end{equation}
Accordingly, using \eqref{Diff-J-V-2}, we infer {\eqref{Approx-Insensitzing-Equivalent} and}
\[
	\left| \left.\dfrac{d }{d \tau}\left({J}\left(\Omega_{\tau\mathbf{V}} \right)\right)\right|_{\tau=0}\right| 
	\leq \varepsilon\norm{\mathbf{V}\cdot \mathbf{n}}_{L^\infty(\partial \Omega_0)}, \forall\,\mathbf{V} \in \mathcal{E},
\]
which concludes the proof of Theorem \ref{Thm-Approx-Case1}, \ma{by using the fact that $W^{1,\infty}(\R^d,\R^d)$ is included into the space of continuous vector fields, so that 
$$
\Vert \mathbf{V}\Vert_{L^\infty(\partial\Omega_0)}\leq \Vert \mathbf{V}\Vert_{\mathscr{C}^0(\overline{\Omega_0})}\leq \Vert \mathbf{V}\Vert_{W^{3,\infty}}(\Omega_0)
$$
by inclusion.}

\subsection{Proof of Theorem \ref{Thm-Case1-Exact-FD}}\label{Subsec-Proof-Exact-Case1}
Before proving Theorem \ref{Thm-Case1-Exact-FD}, let us give an insight of the strategy of our proof. It will be divided into three steps.
\begin{itemize}
\item First step: we will treat the case of the exact {insensitizing} for $J_h$ with respect to a finite dimensional space $\mathcal E$ of dimension $1$, in the sense of \eqref{Finite-Insensitizing}, in order to explain the main idea behind our proof. This case was already studied and analyzed with different techniques in \cite{Lissy-Privat-Simpore}.
\item Second step: we will explain how to modify our first step to the case of the exact {insensitizing} for $J_h$ with respect to any finite dimensional space $\mathcal E$, in the sense of \eqref{Finite-Insensitizing}, by using the Brouwer fixed point Theorem.
\item Last step: we will explain how the construction made in the previous step together with the use of Proposition \ref{Prop-Approx-Cont-Case1} ensures that {one can simultaneously solve the $\varepsilon$-approximate insensitizing of $J_h$ and its exact insensitizing with respect to a finite dimensional space $\mathcal{E}$.}
\end{itemize}
Remind that we assume the geometric setting \eqref{Case-1-emptyset}.
\paragraph{First step: Exact insensitizing in the case $\mathcal{E}=\hbox{Span\,} \{\mathbf{V}\} $.}
Let us fix some $\mathbf{V}\in W^{3,\infty}(\mathbb R^d;\mathbb R^d)$ supposed to be non-zero, and let us consider the case $\mathcal{E}=\hbox{Span\,} \{\mathbf{V}\}$, \textit{i.e.} the case of a one-dimensional vector space $\mathcal{E}$, and only focus on the proof of exact insensitizing of $J_h$ with respect to $\mathcal{E} = \hbox{Span\,} \{\mathbf{V}\}$.

Recall that, according to \eqref{Diff-J-V-2}, the exact insensitizing problem for $\mathcal{E}$ is equivalent to determining a control function $h \in L^2(0,T; L^2(\omega_0))$ such that
\begin{equation}
	\label{argenton}
		\int_{\partial\Omega_{0}}(\mathbf{V}\cdot \mathbf{n})
			\left(\int_{0}^{T}\partial_{n}y_{0}\partial_{n} z_{0}\, dt \right)\,
			d\sigma
			=0, 
\end{equation}
where $(y_0,z_0)$  solves \eqref{y-Omega-0-h}.

Of course, if $\mathbf{V}\cdot \mathbf{n}=0$ on \red{$\partial\Omega_0$} (which may happen since we only assumed that $\mathbf{V} $ is non-zero as a function defined in $\mathbb R^d$), then \eqref{argenton} is automatically verified and the problem is trivial. Hence, from now on, we assume that $\mathbf{V}\cdot \mathbf{n}$ does not vanish {identically} on \red{$\partial\Omega_0$}.

To study condition \eqref{argenton}, we introduce the pairs $(y_\xi,z_\xi)$ and $(y_h,z_h)$ as the solutions of the linear systems
\begin{equation}
	\left\{ 
		\begin{array}{rcll}
			\dfrac{\partial y_\xi}{\partial t}-\Delta y_\xi &=&\xi & \text{ in } (0,T) \times \Omega_0 ,
			\\ 
			y_\xi &=& 0 & \text{ on }(0,T) \times \partial \Omega_0,
			\\ 
			y_\xi\left( 0,\cdot \right) &=& 0 & \text{ in }\Omega_{0},
		\end{array}%
	\right.
	\label{cou1b}
\end{equation}
\begin{equation}
	\left\{ 
		\begin{array}{rcll}
			-\dfrac{\partial z_\xi}{\partial t}-\Delta z_\xi &=& y_\xi\chi_{\Theta}& \text{ in } (0,T) \times \Omega_0 ,
			\\ 
			z_\xi&=& 0 & \text{ on }(0,T) \times \partial \Omega_0,
			\\ 
			z_\xi\left( T,\cdot\right) &=& 0 & \text{ in }\Omega_{0},
		\end{array}%
	\right.
	\label{cou2b}
\end{equation}
and
\begin{equation}
	\left\{ 
		\begin{array}{rcll}
			\dfrac{\partial y_h}{\partial t}-\Delta y_h &=& h\mathbbm{1}_{\omega} & \text{ in } (0,T) \times \Omega_0  ,
			\\ 
			y_h &=& 0 & \text{ on }(0,T) \times \partial \Omega_0,
			\\ 
			y_h \left( 0,\cdot \right) &=& 0 & \text{ in }\Omega_{0},
		\end{array}%
	\right.
	\label{cou1bis}
\end{equation}
\begin{equation}
	\left\{ 
		\begin{array}{rcll}
			-\dfrac{\partial z_h}{\partial t}-\Delta z_h  &=& {y_h}\chi_{\Theta}& \text{ in }(0,T) \times \Omega_0  ,
			\\ 
			z_h &=& 0 & \text{ on }(0,T) \times \partial \Omega_0,
			\\ 
			z_h\left( T,\cdot\right) &=& 0 & \text{ in }\Omega_{0}.
		\end{array}%
	\right.
	\label{cou2bis}
\end{equation}
This allows to decompose the solution $(y_0,z_0)$ of \eqref{y-Omega-0-h} as 
\[
	y_{0}=y_\xi+y_h, \quad\text{ and } \quad z_{0}=z_\xi+z_h.
\]
Now, we introduce the function $\mathcal{U}: L^2(0,T; L^2(\omega)) \to \R$ defined for $h \in L^2(0,T; L^2(\omega))$  by 
\begin{equation*}
	%\label{Def-U}
	\mathcal{U}(h) = 
			\int_{\partial\Omega_{0}}(\mathbf{V}\cdot \mathbf{n})
				\left(\int_{0}^{T}(\partial_{n}y_{\xi} + \partial_n y_h) (\partial_{n} z_{\xi} + \partial_n z_h)\, dt \right)\,
			d\sigma, 
\end{equation*}
 so that condition \eqref{argenton} can be simply reformulated as $\mathcal{U}(h) = 0$.
 
Our goal is to find $h \in L^2(0,T; L^2(\omega))$ such that $\mathcal{U}(h)= 0$. \red{To construct such a function $h$}, we will look for a two-dimensional vector space spanned by two elements $h_1$ and $h_2$ in $L^2(0,T; L^2(\omega))$ such that the function $\mathcal{U}$ vanishes \red{for at least one element of} $\text{Span\,} \{ h_1, h_2\}$, \textit{i.e.} we want to show that 
\begin{equation}
	\label{Reformulation-U-h=0}
	\exists (h_1,h_2)\in L^2(0,T; L^2(\omega))^2,\,\exists (\lambda_1, \lambda_2) \in \R^2 \, \text{ such that } \mathcal{U}(\lambda_1 h_1 + \lambda_2 h_2) = 0.
\end{equation}

To show that this can be done, we observe that the map $h \in L^2(0,T; L^2(\omega)) \to (\partial_n y_h, \partial_n z_h)\in (L^2(0,T; L^2(\partial\Omega_0)))^2$ is linear, hence it is obvious that the function $\mathcal{U}$ can be decomposed as follows
\begin{equation*}
%	\label{Decomposition-U}
	\mathcal{U}(h)=Q(h)+L(h)+C,
\end{equation*}
where $Q$ is quadratic in $h$, $L$ is linear in $h$ and $C$ does not depend on $h$:
\begin{eqnarray}
	Q(h)&=& \int_{\partial \Omega_0} (\mathbf{V}\cdot \mathbf{n})
	\left( \int_0^T \partial_{n} y_{h}\partial_{n} z_{h}\, dt\right) d\sigma,
	\label{quadra}
	\\
	L(h)&=& \int_{\partial \Omega_0}(\mathbf{V}\cdot \mathbf{n}) \
	\left(\int_0^T (\partial_{n} y_{\xi}\partial_{n} z_{h}+\partial_{n} y_{h}\partial_{n} z_{\xi})\, dt \right)d\sigma,
	\notag
	\\
	C&=& \int_{\partial \Omega_0} (\mathbf{V}\cdot \mathbf{n})
	\left(\int_0^T \partial_{n} y_{\xi}\partial_{n} z_{\xi}\, dt \right)
	d\sigma .
	\notag
\end{eqnarray}
Accordingly, problem \eqref{Reformulation-U-h=0} amounts to finding $h_1,\,h_2$ in $L^2(0,T; L^2(\omega))$ and $\lambda_1,\,\lambda_2$ in $\R$ such that 
 \begin{equation}
 	\label{eqM1}
		 \lambda_1^2Q_{11}(h_1)+\lambda_1\lambda_2 Q_{12}(h_1,h_2)+\lambda_2^2Q_{22}(h_2)
		 +\lambda_1 L_{1}(h_1)+\lambda_2 L_2(h_2)+C=0,
 \end{equation}
 where 
 \begin{eqnarray*}
	Q_{11}(h_1)&=& \int_{\partial \Omega_0} (\mathbf{V} \cdot \mathbf{n})
	\left(\int_0^T \partial_{n} y_{h_1}\partial_{n} z_{h_1}\, dt\right) d\sigma,
	\\
	Q_{12}(h_1,h_2)&=& \int_{\partial \Omega_0}(\mathbf{V}\cdot \mathbf{n}) \left(\int_0^T (\partial_{n} y_{h_1}\partial_{n} z_{h_2}+\partial_{n} y_{h_2}\partial_{n} z_{h_1})\, dt\right) d\sigma,
	\\
	Q_{22}(h_2) &=& \int_{\partial \Omega_0} (\mathbf{V} \cdot \mathbf{n})\left(\int_0^T \partial_{n} y_{h_2}\partial_{n} z_{h_2}\, dt\right) d\sigma,
	\\
	L_1(h_1) &=& \int_{\partial \Omega_0} (\mathbf{V} \cdot \mathbf{n})\left(\int_0^T (\partial_{n} y_{h_1}\partial_{n} z_{\xi}+\partial_{n} y_{\xi}\partial_{n} z_{h_1})\, dt\right) d\sigma,
	\\
	L_2(h_2) &=& \int_{\partial\Omega_0} (\mathbf{V} \cdot \mathbf{n})\left(\int_0^T(\partial_{n} y_{h_2}\partial_{n} z_{\xi}+\partial_{n} y_{\xi}\partial_{n} z_{h_2})\, dt \right)d\sigma.
	%\\
	%F &=& \int_{(0,T)\times\partial\Omega_0} (\mathbf{V}_i\cdot \mathbf{n})\partial_{n} y_{\xi}\partial_{n} q_{\xi}\, dt d\sigma.
\end{eqnarray*}

Our strategy then reduces to choose $h_1$ and $h_2$ such that the  Neumann traces $(\partial_n y_{h_i},\partial_n z_{h_i})$, $i=1,2$ in $L^2(0,T; L^2(\partial\Omega_0))$ for the solutions of   \eqref{cou1bis}-\eqref{cou2bis} with $h_i$, allows to guarantee the existence of a solution $(\lambda_1, \lambda_2) \in \R^2$ to \eqref{eqM1}.

Let us choose $(\gamma_{i,y}, \gamma_{i,z})$ in $(L^2(0,T; L^2(\partial \Omega_0)))^2$ for $i = 1, 2$ as follows: 
\begin{align*}
	&\gamma_{1,y} =  \frac{(\mathbf{V} \cdot \mathbf{n})}{\norm{\mathbf{V}\cdot \mathbf{n}}_{L^2(\partial\Omega_0)}^2} \frac{\mathbbm{1}_{(0,T/2)} }{T}, 
	\quad
	&&
	\gamma_{1,z} = \mathbbm{1}_{(T/2,T)}, 
	\\
	&
	\gamma_{2, y} = \frac{(\mathbf{V} \cdot \mathbf{n})}{\norm{\mathbf{V}\cdot \mathbf{n}}_{L^2(\partial\Omega_0)}^2} \frac{\mathbbm{1}_{(T/2,T)} }{T}, 
	\quad
	&&\gamma_{2, z} = \mathbbm{1}_{(0,T/2)}.
\end{align*}

We easily have that 
\begin{align*}
	& \int_0^T \gamma_{1, y} (t, x) \gamma_{1,z}(t,x) \, dt = 0, \hbox{ for all } x \in \partial \Omega_0, 
	\\
	& \int_0^T \gamma_{2, y} (t, x) \gamma_{2,z}(t,x) \, dt = 0, \hbox{ for all } x \in \partial \Omega_0, 
	\\
	&
	\int_0^T(\gamma_{1,y} (t,x) \gamma_{2,z} (t,x)+ \gamma_{2,y} (t,x)\gamma_{1,z} (t,x))\, dt = \frac{(\mathbf{V} \cdot \mathbf{n})(x)}{\norm{\mathbf{V}\cdot \mathbf{n}}_{L^2(\partial\Omega_0)}^2},\, \hbox{ for all } x \in \partial \Omega_0.
\end{align*}

If it was possible to find some $h_1\in L^2(0,T; L^2(\omega))$ and $h_2\in L^2(0,T; L^2(\omega))$ such that 
\begin{equation}
	\label{exh}
	(\partial_{n} y_{h_1},\partial_{n} z_{h_1},\partial_{n} y_{h_2},\partial_{n}z_{h_2})
	=(\gamma_{1,y},\gamma_{1,z},\gamma_{2,y},\gamma_{2,z}),
\end{equation}
then, we would have $Q_{11}=Q_{22}=0$ and $Q_{12}=1$, so that  equation \eqref{eqM1} with $\lambda_2 = |\lambda_1|$ would become
\begin{equation*}
%	\label{eqM12}
	\lambda_1|\lambda_1|+\lambda_1 L_1(h_1)+|\lambda_1| L_2(h_2)+C=0,
\end{equation*}
which can obviously be solved for some $\lambda_1 \in \R$ according to the intermediate value theorem, since the left hand-side goes to $-\infty$ when $\lambda_1 \to -\infty$ and to $+\infty$ when $\lambda_1 \to +\infty$ while being continuous on $\R$.

Unfortunately, we cannot \emph{a priori} find $h_1\in L^2(0,T; L^2(\omega))$ and $h_2\in L^2(0,T; L^2(\omega))$ such that \eqref{exh} \textit{exactly} holds, but Proposition \ref{Prop-Approx-Cont-Case1} ensures that \eqref{exh}  \textit{approximately} holds, in the following sense: for any ${\alpha}>0$, there exists  $h_1^{\alpha} \in L^2(0,T; L^2(\omega))$ and $h_2^{\alpha} \in L^2(0,T; L^2(\omega))$ such that 
\begin{equation*}
%\label{app}
	\begin{aligned}
		\norm{\partial_{n} y_{h_1^{\alpha}}-\gamma_{1,y}}_{L^2(0,T; L^2(\partial \Omega_0))}\leqslant {\alpha}, 
		\\
		\norm{\partial_{n} q_{h_1^{\alpha}}-\gamma_{1,q}}_{L^2(0,T; L^2(\partial \Omega_0))}\leqslant {\alpha}, 
		\\
		\norm{\partial_{n} y_{h_2^{\alpha}}-\gamma_{2,y}}_{L^2(0,T; L^2(\partial \Omega_0))}\leqslant {\alpha},
		\\ 
		\norm{\partial_{n} q_{h_2^{\alpha}}-\gamma_{2,q}}_{L^2(0,T; L^2(\partial \Omega_0))}\leqslant {\alpha}.
	\end{aligned}
\end{equation*}

Accordingly, with this choice of $h_1^{\alpha}$ and $h_2^{\alpha}$, the quadratic part $Q$ that is given in \eqref{quadra} is only slightly perturbed in the sense that  
\[
	|Q_{11}(h_1^{\alpha})| +|Q_{12}(h_1^{\alpha},h_2^{\alpha}) -1| + |Q_{22} (h_2^{\alpha})| \leq C {\alpha},
\]
where $C$ only depends on the norm of $(\gamma_{i,y}, \gamma_{i,z})_{i \in \{1,2\}}$ in $L^2(0,T; L^2(\partial \Omega_0))$. Therefore, taking ${\alpha} >0$ such that $C {\alpha} \leq 1/2$, and choosing $(h_1, h_2 ) = (h_1^{\alpha}, h_2^{\alpha})$, we get that 
\[
	|Q(\lambda_1h_1 +|\lambda_1|h_2) - \lambda_1 |\lambda_1| | \leq \frac{|\lambda_1|^2}{2}.
\]
Accordingly, the continuous function $\lambda_1 \in \R \mapsto Q(\lambda_1h_1 +|\lambda_1|h_2)$ goes to $-\infty$ as $\lambda_1 \to -\infty$ and to $+ \infty$ as $\lambda_1 \to \infty$, and hence, the function
\[
	\lambda_1 \in \R \mapsto \mathcal{U}(\lambda_1 h_1+|\lambda_1| h_2)
\]
inherits the same property. Hence, it vanishes for some $\lambda_1 \in \R$. 

This concludes the proof of exact insensitizing of $J_h$ for a vector space $\mathcal{E}$ is of dimension $1$.

\paragraph{Second step: Exact insensitizing for a finite-dimensional vector space $\mathcal{E}$.} Now, we assume that $\mathcal{E}$ is of finite dimension $N\geqslant 2$. Our goal is to mimic the method developed when $\mathcal{E}$ was a one-dimensional vector space, replacing the intermediate value theorem by a Brouwer fixed point argument. 

Let $\mathscr{E} = \{ \mathbf{V}\cdot \mathbf{n},\, \mathbf{V} \in \mathcal{E}\}$, which is itself a finite dimensional subspace of $L^2(\partial\Omega_0)$ of dimension $M\leqslant N$, and choose an orthonormal basis $(\mathbf{V}_k\cdot \mathbf{n})_{k \in \llbracket 1, M\rrbracket}$ of $\mathscr{E}$ for the canonical inner product on $L^2(\partial\Omega_0)$. 

Following the previous case, for all $k \in \llbracket 1,M \rrbracket $, we introduce 
\begin{equation}
	\label{argenton2}
	\mathcal{U}_k(h):=
		\int_{\partial\Omega_{0}}(\mathbf{V}_k\cdot \mathbf{n})
			\left(\int_{0}^{T}(\partial_{n} y_{\xi} + \partial_n y_h) (\partial_{n} z_{\xi} + \partial_n z_h)\, dt \right)
			\, d\sigma, 
\end{equation}
where $y_{h}$, $z_{h}$, $ y_{\xi}$ and $z_{\xi}$ are defined in \eqref{cou1b}, \eqref{cou2b}, \eqref{cou1bis} and \eqref{cou2bis}.

According to \eqref{Finite-Insensitizing}, the insensitizing problem for $J_h$ for the family $\mathcal{E}$ amounts to finding {a function} $h \in L^2(0,T; L^2(\omega))$ such that for all $k \in \llbracket 1,M \rrbracket $, $\mathcal{U}_k(h) = 0$. 

As in the first step, for all $k \in \llbracket 1,M \rrbracket $, the function $\mathcal{U}_k$ can be decomposed as 
\begin{equation}
	\label{Decomposition-U-k}
	\mathcal{U}_k(h)=Q_k(h)+L_k(h)+C_k,
\end{equation}
where $Q_k$ is quadratic in $h$, $L_k$ is linear in $h$ and $C_k$ does not depend on $h$:
\begin{eqnarray}
	\label{Def-Q-k}
	Q_k(h)&=& \int_{\partial \Omega_0} (\mathbf{V}_k\cdot \mathbf{n})
	\left( \int_0^T \partial_{n} y_{h}\partial_{n} z_{h}\, dt\right) d\sigma,
	\\
	\label{Def-L-k}
	L_k(h)&=& \int_{\partial \Omega_0}(\mathbf{V}_k\cdot \mathbf{n}) \
	\left(\int_0^T (\partial_{n} y_{\xi}\partial_{n} z_{h}+\partial_{n} y_{h}\partial_{n} z_{\xi})\, dt \right)d\sigma,
	\\
	\label{Def-C-k}
	C_k&=& \int_{\partial \Omega_0} (\mathbf{V}_k\cdot \mathbf{n})
	\left(\int_0^T \partial_{n} y_{\xi}\partial_{n} z_{\xi}\, dt \right)
	d\sigma .
\end{eqnarray}

For each $k \in \llbracket 1,M \rrbracket $, we introduce the following elements of $L^2(0;T; L^2(\partial \Omega_0))$:
\begin{align*}
	&\gamma_{k,1,y} =  \frac{(\mathbf{V}_k \cdot \mathbf{n})}{\norm{\mathbf{V}_k\cdot n}_{L^2(\partial\Omega_0)}^2} \frac{M }{T} \mathbbm{1}_{((k-1)T/M,(2k-1)T/(2M))}, 
	\quad
	&& \gamma_{k,1,z} = \mathbbm{1}_{((2k-1)T/(2M), kT/M)}, 
	\\
	&
	\gamma_{k,2, y} = \frac{(\mathbf{V}_k \cdot \mathbf{n})}{\norm{\mathbf{V}_k\cdot n}_{L^2(\partial\Omega_0)}^2} \frac{M}{T} \mathbbm{1}_{((2k-1)T/(2M), kT/M)}, 
	\quad
	&&\gamma_{k,2, z} = \mathbbm{1}_{((k-1)T/M,(2k-1)T/(2M))}.
\end{align*}

It is then easy to check that 
\begin{multline}
	\label{Identity-gamma-i-j-k}
	\forall (i,j,k) \in \llbracket 1,M \rrbracket^3 , \, (a,b) \in \{1,2\}^2, 
	\\
	\int_{\partial \Omega_0} \mathbf{V}_k \cdot \mathbf{n} \left( \int_0^T (\gamma_{i,a, y}  \gamma_{j,b,z}  + \gamma_{j,b,y} \gamma_{i,a,z}) \, dt \right) d \sigma 
	= \delta_{i,j,k} \mathbbm{1}_{a \neq b},
\end{multline}
where $ \delta_{i,j,k}$ denotes {the Kronecker symbol ($\delta_{i,j,k} = 1$ if and only if $i = j = k$, and $=0$ otherwise).}

Now, for $k \in \llbracket 1,M \rrbracket $ and $a \in \{1, 2\}$, using Proposition \ref{Prop-Approx-Cont-Case1}, for any ${\alpha} >0$, there exists $h_{k,a}^{{\alpha}} \in L^2(0,T; L^2(\omega))$ such that the solution $(y_{h_{k,a}^{\alpha}}, z_{h_{k,a}^{\alpha}})$ of \eqref{cou1bis}--\eqref{cou2bis} satisfies: 
\begin{equation}
	\label{Almost-gamma-k-a}
	\norm{\partial_n y_{h_{k,a}^{\alpha}} - \gamma_{k,a,y}}_{L^2(0,T; L^2(\partial \Omega_0))}
	+ 
	\norm{\partial_n z_{h_{k,a}^{\alpha}} - \gamma_{k,a,y}}_{L^2(0,T; L^2(\partial \Omega_0))}
	\leq 
	{\alpha}.
\end{equation}

Using \eqref{Identity-gamma-i-j-k} and \eqref{Almost-gamma-k-a}, we easily show that, for any $\lambda = (\lambda_{k})_{k \in \llbracket 1,M \rrbracket } \in \R^{M}$, for any $k \in \llbracket 1,M \rrbracket $,
\[
	\left| 
		Q_k \left( \sum_{j = 1}^M (\lambda_{j} h_{j,1}^{\alpha} + |\lambda_{j}| h_{j,2}^{\alpha}) \right)
		- 
		\lambda_k |\lambda_k|
	\right|
	\leq 
	C {\alpha} \norm{\lambda}_{\R^M}^2, 
\] 
for some $C >0$ independent of ${\alpha}$, where $Q_k$ is defined in \eqref{Def-Q-k}. Therefore, choosing ${\alpha} >0$ small enough such that $C {\alpha} \leq 1/(2M)$ and dropping the superscript ${\alpha}$ from now on, we have, for all $\lambda = (\lambda_{k})_{k \in \llbracket 1,M \rrbracket } \in \R^{M}$, for any $k \in \llbracket 1,M \rrbracket $,
\begin{equation}
	\label{Key-Estimate-Q-k}
	\left| 
		Q_k  \left( \sum_{j = 1}^M (\lambda_{j} h_{j,1} + |\lambda_{j}| h_{j,2}) \right)
		- 
		\lambda_k |\lambda_k|
	\right|
	\leq 
	\frac{1}{2M}  \norm{\lambda}_{\R^M}^2, 
\end{equation}

Our next goal is to check that 
\begin{equation}
	\label{Goal-M}
	\exists (\lambda_{k})_{k \in \llbracket 1,M \rrbracket } \in \R^{M}
	\text{ such that for all } 
k \in \llbracket 1,M \rrbracket , 
	\quad 
	\mathcal{U}_k\left(  
				 \sum_{j= 1}^M (\lambda_{j} h_{j,1} + |\lambda_{j}| h_{j,2})
				\right)
	= 0, 
\end{equation}
where $\mathcal{U}_k$ is defined in \eqref{argenton2}. Based on the decomposition \eqref{Decomposition-U-k}, in order to do that, we will use a fixed point argument. We introduce the continuous function $s:\mathbb R\rightarrow \mathbb R$ given by 
$$s(y)=\sqrt{y} \mbox{ if }y\geqslant 0,\, s(y)=-\sqrt{-y} \mbox{ if } y<0.$$
Now, let us define the mapping
$$
F : \lambda = (\lambda_{k})_{k \in \llbracket 1,M \rrbracket } \in \R^{M} \longmapsto \widehat \lambda = (\widehat{\lambda_{k}})_{k \in \llbracket 1,M \rrbracket } \in \R^{M} ,
$$
where, for all $k \in \llbracket 1,M \rrbracket $, $\widehat{\lambda_k}$ is defined  as 
$$
\widehat{\lambda_k}= s \left(- \left( Q_k\left( \sum_{j= 1}^M (\lambda_{j} h_{j,1} + |\lambda_{j}| h_{j,2})\right)- \lambda_k |\lambda_k| \right)- L_k \left(\sum_{j= 1}^M (\lambda_{j} h_{j,1} + |\lambda_{j}| h_{j,2})\right)- C_k\right),$$
for $Q_k$, $L_k$ and $C_k$ given in \eqref{Def-Q-k}--\eqref{Def-L-k}--\eqref{Def-C-k}. Notably, by definition of $s$, we have
\begin{equation}\label{chlk}
\widehat{\lambda_k} |\widehat{\lambda_k}| = - \left( Q_k\left( \sum_{j= 1}^M (\lambda_{j} h_{j,1} + |\lambda_{j}| h_{j,2})\right)- \lambda_k |\lambda_k| \right)- L_k \left(\sum_{j= 1}^M (\lambda_{j} h_{j,1} + |\lambda_{j}| h_{j,2})\right)- C_k.
\end{equation}

Using the identity \eqref{chlk}, one easily checks that if $(\lambda_{k})_{k \in \llbracket 1,M \rrbracket } \in \R^{M} $ is a fixed point of $F$, then, it solves the problem \eqref{Goal-M}. 

It is clear that the function $F$ is continuous in $\R^M$. We will simply show that it maps a ball into itself and conclude using Brouwer fixed point theorem. In order to prove that, we use the bound \eqref{Key-Estimate-Q-k} and the following straightforward estimates: for all $k \in \llbracket 1,M \rrbracket $, 
\begin{align*}
	& \left| L_k\left(  
				 \sum_{j= 1}^M (\lambda_{j} h_{j,1} + |\lambda_{j}| h_{j,2})
				\right)
				\right|
				\leq C \norm{(\partial_n y_\xi, \partial_n z_\xi)}_{(L^2(0,T; L^2(\partial \Omega_0)))^2} \norm{\lambda}_{\R^M}, 
			\\
	& | C_k | \leq \tilde{C}\norm{(\partial_n y_\xi, \partial_n z_\xi)}_{(L^2(0,T; L^2(\partial \Omega_0)))^2}^2,  	
\end{align*}
where \textcolor{red}{$\tilde{C}$} is a positive constant (independent of $k \in \llbracket 1,M \rrbracket $ and of $\xi$). 

Accordingly, \red{by denoting $\Vert \cdot \Vert_{\R^M}$ the euclidean norm,} one has for all $\lambda \in \R^M$, 
\begin{align*}
	\norm{\widehat \lambda}_{\R^M}^2
	&\leq 
	\frac{1}{2} \norm{ \lambda}_{\R^M}^2
	+ 
	\tilde{C} \norm{(\partial_n y_\xi, \partial_n z_\xi)}_{(L^2(0,T; L^2(\partial \Omega_0)))^2} \norm{\lambda}_{\R^M}
	+ 
	\tilde{C} \norm{(\partial_n y_\xi, \partial_n z_\xi)}_{(L^2(0,T; L^2(\partial \Omega_0)))^2}^2
	\notag
	\\
	& \leq
	\frac{3}{4} \norm{\lambda}_{\R^M}^2
	+ 
	2\tilde{C}\norm{(\partial_n y_\xi, \partial_n z_\xi)}_{(L^2(0,T; L^2(\partial \Omega_0)))^2}^2.
\end{align*}
In particular, setting 
\[
	R  = 2\sqrt{2\tilde C} \norm{(\partial_n y_\xi, \partial_n z_\xi)}_{(L^2(0,T; L^2(\partial \Omega_0)))^2}, 
\]
\red{where $\tilde C$ denotes the constant in the previous estimate,} the closed ball of $\R^M$ of radius $R$ \red{(for the euclidean norm)} is stable by $F$. Therefore, by Brouwer fixed point theorem, there exists $\lambda \in \R^M$ in the closed ball of radius $R$ such that $F(\lambda) = \lambda$.

This proves the existence of $\lambda = (\lambda_k)_{k \in \llbracket 1,M \rrbracket } \in \R^M$ satisfying \eqref{Goal-M} and with the bound 
\begin{equation}
	\label{Est-Lambda}
	\norm{\lambda}_{\R^M} \leq \widehat C \norm{(\partial_n y_\xi, \partial_n z_\xi)}_{(L^2(0,T; L^2(\partial \Omega_0)))^2}\red{,}
\end{equation}
\red{for some positive constant $\widehat C$.}
In particular, this bound implies that the corresponding control 
\begin{equation}
	\label{Control-h-Exact}
	h =  \sum_{j= 1}^M (\lambda_{j} h_{j,1} + |\lambda_{j}| h_{j,2}), 
\end{equation}
and the corresponding controlled trajectory $(y_h,z_h)$ of \eqref{cou1bis}--\eqref{cou2bis} satisfies
\begin{multline}
	\label{Bound-p-n-y-h-z-h}
	\norm{h}_{L^2(0,T; L^2(\omega))} 
	+ 
	\norm{(y_h,z_h)}_{(L^2(0,T; H^2(\Omega_0)))^2} 
	+ 
	\norm{(\partial_n y_h, \partial_n z_h)}_{(L^2(0,T; L^2(\partial \Omega_0)))^2}
	\\
	\leq 
	C \norm{(\partial_n y_\xi, \partial_n z_\xi)}_{(L^2(0,T; L^2(\partial \Omega_0)))^2},
\end{multline}
for some positive constant $C$.

\paragraph{Last step: adding the approximate insensitizing property.} The idea there is to decompose the approximate and exact insensitizing problems. {To be more precise, we will choose the control $h$ in two steps, under the form $h = h_0 + h_1$, where $h_0$ is used to get the approximate insensitizing property, and $h_1$ is then chosen afterwards to get the exact insensitizing property in the directions of $\mathcal{E}$.}

Given $\xi \in L^2(0,T; L^2(\Omega_0))$ and $\varepsilon >0$, we set 
\begin{equation*}
%	\label{Def-varepsilon-0}
	\varepsilon_0 = \frac{{\sqrt{\varepsilon}}}{C+1}, \text{ where $C$ is the constant in \eqref{Bound-p-n-y-h-z-h}}.
\end{equation*}
According to Proposition \ref{Prop-Approx-Cont-Case1}, there exists $h_0 \in L^2(0,T; L^2(\omega))$ such that the solution $(y,z)$ of \eqref{Coupled-Heat-Approx} satisfies 
\begin{equation}
	\label{Approx-Result-First}
	 \norm{(\partial_n y, \partial_n z)}_{(L^2(0,T; L^2(\partial \Omega_0)))^2}
	 \leq 
	 \varepsilon_0.
\end{equation}
Now, we set $\xi_1 = \xi + h_0 \mathbbm{1}_\omega$, which belongs to $L^2(0,T; L^2(\Omega_0))$. According to the previous paragraph applied for the source term $\xi_1$, there exists $h_1 \in L^2(0,T;L^2(\omega))$ such that 
\begin{equation}
	\label{Exact-insensitizing-E}
	\forall \mathbf{V} \in \mathcal{E}, 
	\quad 
	 \int_{\partial\Omega_{0}}(\mathbf{V}\cdot \mathbf{n})
		\left(\int_{0}^{T}\partial_{n}y_{0}\partial_{n} z_{0}\, dt \right)\,
		d\sigma
	=0, 
\end{equation}
where $(y_0, z_0)$ is the solution of 
\begin{equation}
	\label{y-Omega-00}
		\left\{ 
			\begin{array}{rcll}
				\dfrac{\partial y_0}{\partial t}-\Delta y_0 &=&\xi_1 +h_1\mathbbm{1}_{\omega
}= \xi + (h_0 + h_1) \mathbbm{1}_\omega & \text{ in }(0,T) \times \Omega_0,
				\\ 
				y_0 &=& 0 & \text{ on }(0,T) \times \partial \Omega_0,
				\\ 
				y_0\left( 0,\cdot\right) &=& 0 & \text{ in }\Omega_0
				\\
				-\dfrac{\partial z_0}{\partial t}-\Delta z_0 &=&\mathbbm{1}_{\Theta}y_0 & \text{ in }(0,T) \times \Omega_0,
				\\ 
				z_0 &=& 0 & \text{ on }(0,T) \times \partial \Omega_0,
				\\ 
				z_0\left( T,\cdot\right) &=& 0 & \text{ in }\Omega_0.
			\end{array}
		\right.
\end{equation}
Besides, 
\[
	\partial_n y_0 = \partial_n y + \partial_n y_{h_1}, \quad \partial_n z_0 = \partial_n z + \partial_n z_{h_1}, 
\]
so that the bounds \eqref{Bound-p-n-y-h-z-h} and \eqref{Approx-Result-First} imply
\[
	 \norm{(\partial_n y_0, \partial_n z_0)}_{(L^2(0,T; L^2(\partial \Omega_0)))^2}
	 \leq 
	(C+1) \varepsilon_0, 
\]
where $C$ is the constant in \eqref{Bound-p-n-y-h-z-h}. We then easily get that  
\begin{align*}
	\forall \mathbf{V} \in W^{3, \infty}(\R^d, \R^d), 
	\quad 
	\left| \int_{\partial\Omega_{0}}(\mathbf{V}\cdot \mathbf{n})
		\left(\int_{0}^{T}\partial_{n}y_{0}\partial_{n} z_{0}\, dt \right)\,
		d\sigma
		\right|
	&\leq (C+1)^2 \varepsilon_0^2  \| \mathbf{V}\cdot \mathbf{n}\|_{L^\infty(\partial \Omega_0)}
	\\
	&\leq \varepsilon  \| \mathbf{V}\cdot \mathbf{n}\|_{L^\infty(\partial \Omega_0)}.
\end{align*}
In other words, $h = h_0 + h_1$ exactly insensitizes $J_h$ for $\mathcal{E}$ and $\varepsilon$-approximately insensitizes $J_h$.

%%%%%%%%%%%%%%%%%%%%%%%%%%%%%%%%%%%%%%%%%%%%%%%%%%%%%%%%%
%
%
%

\section{The case $\omega \cap \Theta \neq \emptyset$}
\label{Sec-Omega-Theta-Not-Empty}

In this whole section, we assume the geometric setting \eqref{Case-2-non-emptyset}.%{, i.e. $\omega \cap \Theta \neq \emptyset$.}

\subsection{Proof of Proposition \ref{Prop-Approx-Cont-Case2}}\label{Subsec-Proof-ApproxCont-Case2}

Similarly to the proof of Proposition \ref{Prop-Approx-Cont-Case1}, we reformulate the first part of Proposition \ref{Prop-Approx-Cont-Case2}, \textit{i.e.} the approximate controllability property \eqref{Approx-Cont-Case2},  as the density of the range of the operator $\mathcal{L}$ defined for $h \in L^2(0,T; L^2(\omega))$ with values in $(L^2(0,T; L^2(\partial \Omega_0)))^2 \times L^2(\Omega_0)$ by 
\[
	\mathcal{L} (h) = ( \partial_n y, \partial_n z, y(T)), 
\]
where $(y,z)$ is the solution of \eqref{csanss}. Arguments similar to the one developed in the proof of Proposition  \ref{Prop-Approx-Cont-Case1} imply that $\mathcal L$ is a bounded operator.
\red{By using similar classical arguments resting upon duality, as in the proof of Proposition~\ref{Prop-Approx-Cont-Case1}, one has $\overline{\text{Ran\,} \mathcal{L} } = (L^2(0,T; L^2(\partial \Omega_0)))^2 \times L^2(\Omega_0)$ if and only if $\text{Ker\,} \mathcal{L}^*  = \{0\}$ \red{(we refer for instance to \cite[Corollary~2.18]{Brezis})}. Using the arguments developed in \cite[Proposition 2.1]{Lissy-Privat-Simpore} together with \cite[Theorem 2.43, page 56]{CNL} (in order to take into account the term $y(T)$ in the definition of $\mathcal L$), one can compute explicitly $\mathcal L^*$ and deduce that $\text{Ker}\,\mathcal{L}^* = \{0\}$ is equivalent to the following unique continuation problem:} if $(g_1, g_2, \psi_T ) \in (L^2(0,T; L^2(\partial\Omega_0)))^2\times L^2(\Omega_0)$, and $(\psi, \varphi)$ solves 
\begin{equation}
	\label{Coupled-Heat-Approx-Adj-Case2}
		\left\{ 
			\begin{array}{rcll}
				- \dfrac{\partial \psi}{\partial t}-\Delta \psi &=&\mathbbm{1}_{\Theta}\varphi & \text{ in }(0,T) \times \Omega_0,
				\\ 
				\psi &=& g_1 & \text{ on }(0,T) \times \partial \Omega_0,
				\\ 
				\psi \left( T,\cdot\right) &=& \psi_T & \text{ in }\Omega_0, 
				\\			
				 \dfrac{\partial \varphi}{\partial t}-\Delta \varphi &=& 0 & \text{ in }(0,T) \times \Omega_0,
				\\ 
				\varphi &=& g_2 & \text{ on }(0,T) \times \partial \Omega_0,
				\\ 
				\varphi\left( 0,\cdot\right) &=&0 & \text{ in }\Omega_0,
			\end{array}
		\right.
\end{equation}
then, we have the following unique continuation property:
\begin{equation}
	\label{UC-Adj-Coupled-Case2}
	\psi = 0 \text{ in } (0,T) \times \omega \Rightarrow g_1 =g_2 =  0 \text{ and } \psi_T = 0.
\end{equation} 

Let us then take $(g_1, g_2, \psi_T ) \in (L^2(0,T; L^2(\partial\Omega_0)))^2\times L^2(\Omega_0)$, $(\psi, \varphi)$ solving \eqref{Coupled-Heat-Approx-Adj-Case2} with $\psi = 0$ in $(0,T) \times \omega$.  Then the equation \eqref{Coupled-Heat-Approx-Adj-Case2} on $\psi$ implies that $\varphi = 0$ in $(0,T)\times {(\omega\cap \Theta)}$. The classical unique continuation property for the heat equation \red{(for instance, by the Holmgren uniqueness theorem, see \textit{e.g.} \cite[Section~8.6]{Hormander-2003})} then applies to $\varphi$ and $\varphi = 0$ in $(0,T)\times \Omega_0$. Following arguments similar to the proof of Proposition Proposition~\ref{Prop-Approx-Cont-Case1}, we deduce that $g_2 = 0$, and that $\psi$ solves the heat equation 
\[
	-\dfrac{\partial \psi}{\partial t}-\Delta \psi =0, \quad \text{ in }(0,T) \times \Omega_0.
\]
Since $\psi = 0$ in $(0,T) \times \omega$, we immediately deduce  \red{again, from the Holmgren uniqueness theorem (see \textit{e.g.} \cite[Section~8.6]{Hormander-2003})}, that $\psi = 0$ in $(0,T) \times \Omega_0$. Since $\psi\in C^0([0,T],H^{-1}(\Omega_0))$ (see \textit{e.g.} \cite[Proposition 7.1.3]{TWBook}), the initial datum $\psi_T$ vanishes as well.  Following arguments similar to the proof of Proposition Proposition~\ref{Prop-Approx-Cont-Case1}, we deduce that $g_1 = 0$.
 
This finishes the proof of the unique continuation property \eqref{UC-Adj-Coupled-Case2} for solutions of \eqref{Coupled-Heat-Approx-Adj-Case2}, hence the proof of the first part of Proposition \ref{Prop-Approx-Cont-Case2}, \textit{i.e.} of the approximate controllability property \eqref{Approx-Cont-Case2}.
\medskip

Let us now focus on the proof of the null-controllability property when $\xi \in L^2(0,T; L^2(\Omega_0))$ can be steered to $0$ using a control $h_{nc} \in L^2(0,T; L^2(\omega))$, in the sense that the solution $y_{nc}$ of \eqref{y-nc} satisfies \eqref{NC-Assumption}. 

Then, for any $\varepsilon >0$ and $f_1$, $f_2$ in $L^2(0,T; L^2(\partial \Omega_0))$, we look for $y$, $z$ and $h$ respectively as $y = y_{nc} + y_1$, $z = z_1$ and $h = h_{nc} + h_1$ where $(y_1, z_1)$ satisfies 
\begin{equation}
	\label{y-Omega-1-h}
		\left\{ 
			\begin{array}{rcll}
				\dfrac{\partial y_1}{\partial t}-\Delta y_1 &=&h_1\mathbbm{1}_{\omega
}& \text{ in }(0,T) \times \Omega_0,
				\\ 
				y_1 &=& 0 & \text{ on }(0,T) \times \partial \Omega_0,
				\\ 
				y_1\left( 0,\cdot\right) &=& 0 & \text{ in }\Omega_0,
				\\
				-\dfrac{\partial z_1}{\partial t}-\Delta z_1 &=&\mathbbm{1}_{\Theta} (y_{nc} + y_1) & \text{ in }(0,T) \times \Omega_0,
				\\ 
				z_1 &=& 0 & \text{ on }(0,T) \times \partial \Omega_0,
				\\ 
				z_1\left( T,\cdot\right) &=& 0 & \text{ in }\Omega_0, 
			\end{array}
		\right.
\end{equation}
and
\begin{equation}
	\| \partial_n y_1 - \tilde f_1 \|_{L^2(0,T; L^2(\partial\Omega_0))} 
	+ 
	\| \partial_n z_1 -  f_2 \|_{L^2(0,T; L^2(\partial\Omega_0))} 
	\leq
	\varepsilon
	\quad \text{ and } \quad 
	y_1(T) = 0 \text{  in } \Omega_0, 
	\label{Goal-NC-and-Approx-Case-2}
\end{equation}
where 
\begin{equation}\label{tilf}\tilde f_1= f_1 - \partial_n y_{nc}.\end{equation}

In order to do that, we consider the functional $K_\varepsilon$ defined for $(g_1, g_2, \psi_T) \in (L^2(0,T; L^2(\red{\partial \Omega_0})))^2 \times L^2(\Omega_0)$ by 
\begin{multline*}
	K_{\varepsilon}(g_1, g_2, \psi_T) 
	= 
	\frac{1}{2} \int_0^T \int_\omega |\psi(t,x)|^2 \, dx dt 
	\\
	+ 
	\int_0^T \int_{\partial \Omega_0} \left( \tilde f_1 g_1 +  f_2 g_2 \right) \, d\sigma dt
	+
	\int_0^T \int_{\Omega_0} \mathbbm{1}_{\Theta} y_0 \varphi
	+
	\varepsilon \norm{(g_1,g_2)}_{(L^2(0,T; L^2(\red{\partial \Omega_0})))^2}, 
\end{multline*}
where $(\psi, \varphi)$ is the solution of \eqref{Coupled-Heat-Approx-Adj-Case2} {corresponding to $(g_1, g_2, \psi_T)$.}

Then, we endow the set 
\[
	X_{0} = (L^2(0,T; L^2(\red{\partial \Omega_0})))^2 \times L^2(\Omega_0)
\]
with the norm
\[
	\norm{(g_1, g_2, \psi_T)}_{obs}^2 
	=
	 \int_0^T \int_\omega |\psi(t,x)|^2 \, dx dt 
	 + 
	 \int_0^T \int_{\partial \Omega_0} \left( | g_1|^2 + |g_2|^2 \right) \, d\sigma dt, 
\]
where $\psi$ is the solution of \eqref{Coupled-Heat-Approx-Adj-Case2}. The fact that this defines a norm comes from the unique continuation property \eqref{UC-Adj-Coupled-Case2}. Then, we define%
\begin{equation*}
	X_{obs} = \overline{X_{0}}^{\norm{\cdot}_{obs}}.
\end{equation*}
We claim that $K_\varepsilon$ can be extended continuously to $X_{obs}$. Indeed, let us emphasize that for $(g_1, g_2, \psi_T) \in X_{obs}$, $\psi|_{(0,T) \times \omega}$ is well defined by density as an element of $L^2(0,T;L^2(\omega))$ since the function $(g_1, g_2, \psi_T) \in X_{0} \mapsto \psi|_{(0,T) \times \omega}$ is well-defined on $X_0$ and continuous (by construction) for the topology of $X_{obs}$, and that we have the following straightforward estimate of the solution $\varphi$ of \eqref{Coupled-Heat-Approx-Adj-Case2}: there exists $C>0$ such that for any $(g_1, g_2, \psi_T) \in X_{0}$, 
\begin{equation}
	\label{Est-Varphi-Direct}
	\norm{ \varphi }_{L^2(0,T; L^2(\Omega_0))} \leq C \norm{ g_2}_{L^2(0,T; L^2(\partial \Omega_0))} \leq C \norm{(g_1, g_2, \psi_T)}_{obs}.
\end{equation}
\red{Now, given $(g_1, g_2, \psi_T) \in X_0$, $\psi$ satisfies the backward heat equation \eqref{Coupled-Heat-Approx-Adj-Case2}. We then introduce the solution $\psi_1$ of 
$$
		\left\{ 
			\begin{array}{rcll}
				- \dfrac{\partial \psi_1}{\partial t}-\Delta \psi_1 &=& 0 & \text{ in }(0,T) \times \Omega_0,
				\\ 
				\psi_1 &=& g_1 & \text{ on }(0,T) \times \partial \Omega_0,
				\\ 
				\psi_1 \left( T,\cdot\right) &=& 0 & \text{ in }\Omega_0, 
			\end{array}
		\right.
$$
which belongs to $L^2(0,T; L^2(\Omega_0))$ with 
$$
	\| \psi_1 \|_{L^2(0,T; L^2(\Omega_0))} \leq C \| g_1\|_{ L^2(0,T; L^2(\partial \Omega_0))}.
$$
Next, $\psi - \psi_1$ satisfies the heat equation with homogeneous Dirichlet boundary conditions, so that classical Carleman estimates, see for instance \cite[Theorem 9.4.1]{TWBook} (after having bounded the weight function $\alpha(x)$ by some constant  $C>0$ from above and by $1$ from below) implies
$$
	\| e^{-C/(t(T-t))} (\psi - \psi_1) \|_{L^2(0,T; L^2(\Omega_0))}
	\leq 
	C \norm{\psi- \psi_1}_{L^2(0,T; L^2(\omega))} 
	+ 
	{C\norm{ - \partial_t \psi - \Delta \psi }_{L^2(0,T; L^2(\Omega_0))}}.
$$
Using that the backward heat equation is well posed, standard energy estimates show that 
$$
	\| e^{-C/(T-t)} (\psi - \psi_1) \|_{L^2(0,T; L^2(\Omega_0))}
	\leq 
	C \norm{\psi- \psi_1}_{L^2(0,T; L^2(\omega))} 
	+ 
	{C\norm{ - \partial_t \psi - \Delta \psi }_{L^2(0,T; L^2(\Omega_0))}}. 
$$
In particular, 
\begin{align*}
%	\label{Estimate-psi}
	\| e^{-C/(T-t)} \psi \|_{L^2(0,T; L^2(\Omega_0))}
	& \leq
	C \norm{\psi}_{L^2(0,T; L^2(\omega))} 
	+ 
	{C\norm{ \varphi }_{L^2(0,T; L^2(\Omega_0))}}
	+ 
	C \norm{g_1}_{L^2(0,T; L^2(\partial \Omega_0))}
	\\
	&\leq 
	C \norm{(g_1, g_2, \psi_T)}_{obs}. 
\end{align*}
Hence, by density of $X_0$ in $X_{obs}$, to each $(g_1, g_2, \psi_T) \in X_{obs}$, we can associate a solution $(\psi, \varphi)$ to
\begin{equation*}
%	\label{Coupled-Heat-Approx-Adj-Case2-Bis}
		\left\{ 
			\begin{array}{rcll}
				- \dfrac{\partial \psi}{\partial t}-\Delta \psi &=&\mathbbm{1}_{\Theta}\varphi & \text{ in }(0,T) \times \Omega_0,
				\\ 
				\psi &=& g_1 & \text{ on }(0,T) \times \partial \Omega_0,
				\\			
				 \dfrac{\partial \varphi}{\partial t}-\Delta \varphi &=& 0 & \text{ in }(0,T) \times \Omega_0,
				\\ 
				\varphi &=& g_2 & \text{ on }(0,T) \times \partial \Omega_0,
				\\ 
				\varphi\left( 0,\cdot\right) &=&0 & \text{ in }\Omega_0,
			\end{array}
		\right.
\end{equation*}
such that for all $T' < T$, $\psi \in L^2(0,T'; L^2(\Omega_0))$. We can thus apply the unique continuation property \eqref{UC-Adj-Coupled-Case2} with $T$ replaced by $T'$, using the same strategy as before relying on Holmgren's uniqueness theorem, and we obtain that if, for $(g_1, g_2, \psi_T) \in X_{obs}$ we have $\psi|_{(0,T') \times \omega} = 0$ then $\varphi$ and $\psi$ vanishes identically on $(0,T')$. Since $T' \in (0,T)$ is arbitrary, we deduce that the unique continuation property \eqref{UC-Adj-Coupled-Case2} extends to $(g_1, g_2, \psi_T) \in X_{obs}$.
}

Classical contradiction arguments relying on \eqref{UC-Adj-Coupled-Case2} (see \textit{e.g.} \cite{Ervedoza-Note-2020}) then give that 
\begin{equation}
	\label{Coercivity-J}
	\liminf_{\norm{(g_1, g_2, \psi_T)}_{obs} \to \infty} \frac{K_\varepsilon(g_1, g_2, \psi_T) }{\norm{(g_1, g_2, \psi_T)}_{obs}} \geq \varepsilon.
\end{equation}

Indeed, let us prove \eqref{Coercivity-J} by contradiction. Assume that there exists a sequence $(g_{1,n}, g_{2,n}, \psi_{T,n}) \in X_{obs}$ such that 
\begin{equation}
	\label{AbsurdumAssumption}
	\rho_n := \norm{(g_{1,n}, g_{2,n}, \psi_{T,n})}_{obs} \to \infty \text{ as } n \to \infty, 
	\text{ and } 
	\alpha:=\limsup_{n \to \infty} \frac{K_\varepsilon(g_{1,n}, g_{2,n}, \psi_{T,n})}{\norm{(g_{1,n}, g_{2,n}, \psi_{T,n})}_{obs}} < \varepsilon.
\end{equation}
We first renormalize the corresponding functions, and set 
\[
	(\tilde g_{1,n}, \tilde g_{2,n}, \tilde \psi_{T,n}) = \frac{(g_{1,n}, g_{2,n}, \psi_{T,n})}{ \norm{(g_{1,n}, g_{2,n}, \psi_{T,n})}_{obs}},
\]
so that 
\begin{equation}
	\label{Bound-tilde-n}
	\norm{(\tilde g_{1,n}, \tilde g_{2,n}, \tilde \psi_{T,n}) }_{obs} = 1.
\end{equation}
Therefore, there exists $(\tilde g_1, \tilde g_2, \tilde \psi_T)$ in $X_{obs}$ such that 
\begin{equation}
	\label{Weak-Convergences}
	(\tilde g_{1,n}, \tilde g_{2,n}, \tilde \psi_{T,n}) \text{ weakly converges to }(\tilde g_1, \tilde g_2, \tilde \psi_T) \text{ in } X_{obs} \text{ as } n \to \infty.
\end{equation}
According to \eqref{AbsurdumAssumption}, we have 
\begin{multline}
	\label{J-eps-expanded}
	\frac{1}{2}(\rho_n)^2 \int_0^T \int_\omega |\tilde \psi_n|^2 \, dx dt 
	\\
	+ 
	\rho_n 
		\left( 
		\int_0^T \int_{\partial \Omega_0} \left( \tilde f_1 \tilde g_{1,n} +  f_2 \tilde g_{2,n} \right) \, d\sigma dt
		+
		\int_0^T \int_{\Omega_0} \chi_{\Theta} y_0 \tilde \varphi_n
		+
		\varepsilon { \norm{(\tilde g_{1,n}, \tilde g_{2,n})}_{(L^2(0,T; L^2(\red{\partial \Omega_0})))^2}
		}
		\right)
	\leq \alpha \rho_n.
\end{multline}
By \eqref{Est-Varphi-Direct} and \eqref{Weak-Convergences}, the sequence 
\[
	\left( 
		\int_0^T \int_{\partial \Omega_0} \left( \tilde f_1 \tilde g_{1,n} +  f_2 \tilde g_{2,n} \right) \, d\sigma dt
		+
		\int_0^T \int_{\Omega_0} \chi_{\Theta} y_0 \tilde \varphi_n
		+
		\varepsilon
		 { \norm{(\tilde g_{1,n}, \tilde g_{2,n})}_{(L^2(0,T; L^2(\red{\partial \Omega_0})))^2}}
	\right)_{n \in \N}
\]
is uniformly bounded in $n$. Therefore, dividing by $\rho_n$ and using that $\rho_n\rightarrow +\infty$ when $n\rightarrow \infty$, we infer
\begin{equation*}
	\lim_{n \to \infty} \left( \int_0^T \int_\omega |\tilde \psi_n|^2 \, dx dt  \right) = 0, 
\end{equation*}
hence 
\[
	\tilde \psi = 0 \hbox{  in } (0,T) \times \omega.
\]
{We also deduce, according to \eqref{Bound-tilde-n}, that 
\begin{equation*}
	\lim_{n \to \infty} \norm{(\tilde g_{1,n}, \tilde g_{2,n})}_{(L^2(0,T; L^2(\red{\partial \Omega_0})))^2} = 1.
\end{equation*}}
{Furthermore, using} the unique continuation property \eqref{UC-Adj-Coupled-Case2}, which has been shown to be also valid for elements of $X_{obs}$, we deduce that 
\[
	(\tilde g_1, \tilde g_2, \tilde \psi_T) = (0,0,0). 
\]
In particular, according to \eqref{Weak-Convergences} and the equation satisfied by $\tilde \varphi$, 
\begin{equation*}
	\lim_{n \to \infty} 
	\left( 
		\int_0^T \int_{\partial \Omega_0} \left( \tilde f_1 \tilde g_{1,n} +  f_2 \tilde g_{2,n} \right) \, d\sigma dt
		+
		\int_0^T \int_{\Omega_0} \chi_{\Theta} y_0 \tilde \varphi_n
		+
		\varepsilon
		 { \norm{(\tilde g_{1,n}, \tilde g_{2,n})}_{(L^2(0,T; L^2(\red{\partial \Omega_0})))^2}}
	\right)
	= 
	\varepsilon,
\end{equation*}
which is in contradiction with \eqref{J-eps-expanded}{, since $\alpha$ was assumed to be smaller than $\varepsilon$ by \eqref{AbsurdumAssumption}}. This concludes the proof of the coercivity estimate \eqref{Coercivity-J}.

Therefore, since the functional $K_\varepsilon$ is also obviously strictly convex \red{on $X_{obs}$}, it admits a unique minimizer $(g_1^*, g_2^*, \psi_T^*) \in X_{obs}$. Writing the Euler-Lagrange equation satisfied by the minimizer gives that, setting 
\begin{equation*}
	h_1 = \psi^* \text{ in } (0,T) \times \omega, 
\end{equation*}
the corresponding solution $(y_1,z_1)$ of \eqref{y-Omega-1-h} satisfies \eqref{Goal-NC-and-Approx-Case-2}. We refer for instance to \cite[Theorem 1.2 and its proof]{Ervedoza-Note-2020} for further details about the approach considered here. To conclude the proof of  Proposition \ref{Prop-Approx-Cont-Case2} , we set $(y,z)=(y_1+y_{nc},z_1)$. Since $(y_1,z_1)$ satisfies \eqref{y-Omega-1-h} and $y_{nc}$ satisfies \eqref{y-nc}, $(y,z)$ is a solution of \eqref{Coupled-Heat-Approx}, verifying moreover \eqref{Approx-Cont-2} and \eqref{NullControl},  thanks to \eqref{Goal-NC-and-Approx-Case-2}, \eqref{tilf} and \eqref{NC-Assumption}.

\subsection{Proof of Theorem \ref{Thm-Approx-Case2}}\label{Subsec-Proof-Approx-Case2}

The proof of Theorem \ref{Thm-Approx-Case2} mainly reduces to Proposition \ref{Prop-Approx-Cont-Case2}, similarly as for the proof of Theorem \ref{Thm-Approx-Case1}, that follows from Proposition \ref{Prop-Approx-Cont-Case1}. 

When the goal is to get the approximate controllability for $y_0$ at time $T>0$, \textit{i.e.} \eqref{Approx-Cont-T} the only novelty is that for any $\varepsilon >0$ and $y_T\in L^2(\Omega_0)$, we should take the control function $h \in L^2(0,T; L^2(\omega))$ such that the solution $(y_0,z_0)$ of \eqref{y-Omega-0-h} satisfies \eqref{Small-Neumann-Data} and \eqref{Approx-Cont-T}, which can be done according to Proposition~\ref{Prop-Approx-Cont-Case2}.

When the goal is to get null-controllability of $y_0$ at time $T$ when $\xi$ is null-controllable, the argument is basically the same, by using a control function $h \in L^2(0,T; L^2(\omega))$ such that the solution $(y_0,z_0)$ of \eqref{y-Omega-0-h} satisfies the estimate \eqref{Small-Neumann-Data}, \textit{i.e.}
$$
	\norm{\partial_n y_0}_{L^2(0,T;L^2(\partial\Omega_0))} + \norm{\partial_n z_0}_{L^2(0,T;L^2(\partial\Omega_0))}\leq \sqrt{\varepsilon}
	{,}
$$
and \eqref{NullControl}. 

Details are left to the reader.

\subsection{Proof of Theorem \ref{Thm-Case2-Exact-FD}}%\label{Subsec-Proof-Exact-FD-Case2}

Here again, the proof of Theorem \ref{Thm-Case2-Exact-FD} is very similar to the one of Theorem \ref{Thm-Case1-Exact-FD}. We focus on the case in which we only want approximate controllability of $y_0$ at time $T$ since the other case in which we want null-controllability at time $T$ when $\xi$ is null-controllable can be deduced similarly, following the proof of Theorem \ref{Thm-Case1-Exact-FD}.

{To be more precise, we will choose the control $h$ in two steps, under the form $h = h_0 + h_1$, where $h_0$ is used to get the approximate insensitizing property and the approximate controllability of $y_0$ at time $T$, and $h_1$ is then chosen afterwards to get the exact insensitizing property in the directions of $\mathcal{E}$.}

\medskip

We assume that $\mathcal{E}$ is of finite dimension $N\geqslant 2$.
Let $\mathscr{E} = \{ \mathbf{V}\cdot \mathbf{n},\, \mathbf{V} \in \mathcal{E}\}$, which is itself a finite dimensional subspace of $L^2(\partial\Omega_0)$ of dimension $M\leqslant N$. Now, we proceed exactly as in the proof of Theorem \ref{Thm-Case1-Exact-FD} by choosing, {for a parameter $\alpha>0$ to be chosen later, }for $k \in \llbracket 1,M \rrbracket $, $a \in \{1,2\}$, control functions $h_{k,a}^\alpha$ such that \eqref{Almost-gamma-k-a} holds and
\begin{equation}\label{yTn1}
	\norm{y_{h_{k,a}^\alpha}(T)}_{L^2(\Omega_0)} \leq 1{,}
\end{equation} 
{where $y_{h_{k,a}^\alpha}$ solves \eqref{cou1bis}--\eqref{cou2bis}. Note that this can be done according to Proposition \ref{Prop-Approx-Cont-Case2}.}

Then, the same arguments as before yields the following result: there exists {$h$ of the form \eqref{Control-h-Exact} with $\lambda \in \R^M$ satisfying \eqref{Est-Lambda}} such that for all $\mathbf{V} \in \mathcal{E}$, 
\[
	\int_{\partial \Omega_0} \mathbf{V} \cdot \mathbf{n} \left( \int_0^T (\partial_n y_\xi + \partial_n y_h) ( \partial_n z_\xi + \partial_n z_h) \, dt \right)\, d\sigma = 0, 
\]
where $(y_\xi, z_\xi)$ solves \eqref{cou1b}--\eqref{cou2b} and $(y_h,z_h)$ solves \eqref{cou1bis}--\eqref{cou2bis}. Besides, combining \eqref{Est-Lambda}, \eqref{Control-h-Exact} and \eqref{yTn1}, the corresponding controlled trajectory $(y_h,z_h)$ of \eqref{cou1bis}--\eqref{cou2bis} satisfies: 
\begin{multline}
	\label{Bound-p-n-y-h-z-h-Case2}
	\norm{h}_{L^2(0,T; L^2(\omega))} 
	+ 
	\norm{(y_h,z_h)}_{(L^2(0,T; H^2(\Omega_0)))^2} 
	+ 
	\norm{y_h(T)}_{L^2(\Omega_0)}
	+ 
	\norm{(\partial_n y_h, \partial_n z_h)}_{(L^2(0,T; L^2(\partial \Omega_0)))^2}
	\\
	\leq 
	C \norm{(\partial_n y_\xi, \partial_n z_\xi)}_{(L^2(0,T; L^2(\partial \Omega_0)))^2}.
\end{multline}
\medskip

Therefore, to solve the problem of approximate insensitizing of $J_h$, exact insensitizing of $J_h$ for $\mathcal{E}$ and approximate controllability \eqref{Approx-Cont-T}, we do as in the proof of Theorem \ref{Thm-Case2-Exact-FD}:  for any $\varepsilon>0$, setting 
\[
	\varepsilon_0 = {\frac{\min\{\sqrt{\varepsilon}, \varepsilon\}}{C + 1},}
	%\sout{\min\left\{ \sqrt{\frac{\varepsilon}{C+1}}, \frac{\varepsilon}{C+1} \right\}}, 
\]
where $C$ is the constant in \eqref{Bound-p-n-y-h-z-h-Case2}, 
we start by taking $h_0 \in L^2(0,T; L^2(\omega))$ such that the solution $(y,z)$ of \eqref{Coupled-Heat-Approx} satisfies %\se{\sout{\eqref{Approx-Result-First} and }}
\[
		\norm{(\partial_n y, \partial_n z)}_{(L^2(0,T; L^2(\partial \Omega_0)))^2}
		+ \norm{y(T) - y_T}_{L^2(\Omega_0)} 
		\leq \varepsilon_0, 
\]
which can be done according to Proposition \ref{Prop-Approx-Cont-Case2}.

Setting $\xi_1 = \xi + h_0 \mathbbm{1}_\omega$, which belongs to $L^2(0,T; L^2(\Omega_0))$, by the previous paragraph and the estimate \eqref{Bound-p-n-y-h-z-h-Case2} applied for the source term $\xi_1$, there exists $h_1 \in L^2(0,T;L^2(\omega))$ such that the identity \eqref{Exact-insensitizing-E} holds for all $\mathbf{V} \in \mathcal{E}$, where $(y_0,z_0)$ denotes the solution of \eqref{y-Omega-00} and the solution $(y_{h_1}, z_{h_1})$ of \eqref{cou1bis}--\eqref{cou2bis} satisfies the bound
\[
	\norm{y_{h_1}(T)}_{L^2(\Omega_0)}
	+ 
	\norm{(\partial_n y_{h_1}, \partial_n z_{h_1})}_{(L^2(0,T; L^2(\partial \Omega_0)))^2}
	\leq 
	C \varepsilon_0.
\]

Besides, 
\[
	y_0(T) = y(T) + y_{h_1}(T) , 
	\quad 
	\partial_n y_0 = \partial_n y + \partial_n y_{h_1},
	\quad 
	\partial_n z_0 = \partial_n z + \partial_n z_{h_1}, 
\]
so that 
\[
	\norm{y_0(T) - y_T}_{L^2(\Omega_0)} +  \norm{(\partial_n y_0, \partial_n z_0)}_{(L^2(0,T; L^2(\partial \Omega_0)))^2}
	 \leq 
	(C+1) \varepsilon_0, 
\]
where $C$ is the constant in \eqref{Bound-p-n-y-h-z-h-Case2}. We then easily get that  
\begin{align*}
	\forall \mathbf{V} \in W^{3, \infty}(\R^d, \R^d), 
	\quad 
	\left| \int_{\partial\Omega_{0}}(\mathbf{V}\cdot \mathbf{n})
		\left(\int_{0}^{T}\partial_{n}y_{0}\partial_{n} z_{0}\, dt \right)\,
		d\sigma
		\right|
	&\leq (C+1)^2 \varepsilon_0^2  \| \mathbf{V}\cdot \mathbf{n}\|_{L^\infty(\partial \Omega_0)}, 
	\\
	& \leq \varepsilon  \| \mathbf{V}\cdot \mathbf{n}\|_{L^\infty(\partial \Omega_0)}, 
\end{align*}
while 
\[
	\norm{y_0(T) - y_T}_{L^2(\Omega_0)} \leq (C+1) \varepsilon_0 \leq \varepsilon.
\]
In other words, $h = h_0 + h_1$ exactly insensitizes $J_h$ for $\mathcal{E}$, $\varepsilon$-approximately insensitizes $J_h$ and $\varepsilon$-approximately controls $y_0$ at time $T$.

%%%%%%%%%%%%%%%%%%%%%%%%%%%%
\section{The exact insensitizing problem}\label{Sec-Exact-insensitizing}

\subsection{Proof of Proposition \ref{Prop-Cas-Facile}}
\label{Sec-Proof-Prop-Cas-Facile}
%%%%%%%%%%%%%%%%%%%%%%%%%%%%
We fix $\xi \in L^2(0,T; L^2(\R^d))$, and we introduce the solution $y_\xi$ of 
\begin{equation*}
%	\label{Heat-Eq-y-f-CasFacile}
	\left\{
		\begin{array}{ll}
			\partial_t y_\xi - \Delta y_\xi = \xi, \quad & (t,x) \in (0,T) \times \Omega_0, 
			\\
			y_\xi(t, x) = 0, \quad & (t,x) \in (0,T) \times \partial \Omega_0, 
			\\
			y_\xi(0,x) = 0, \quad & x \in \Omega_0, 
		\end{array}
	\right.
\end{equation*}
According to \eqref{Assumption-EasyCase}, there exists a smooth function $\eta = \eta(x)$ such that $\eta = 1$ in $\Omega_0 \setminus \omega$ and $\eta=0$ in $\Theta$. Then, we set 
\[
	y_0(t,x) = \eta(x) y_\xi(t,x), \qquad \hbox{ for } (t,x) \in (0,T) \times \Omega_0, 
\]
which solves \eqref{y-Omega-0-h} with control function 
\[
	\mathbbm{1}_{\omega}(x)h(t,x) = (\eta(x) -1) \xi(t,x) - [\Delta, \eta] y_\xi(t,x), \qquad \hbox{ for } (t,x) \in (0,T) \times \Omega_0, 
\]
where $[\Delta, \eta] y_\xi(t,x):=\Delta (\eta y_\xi)-\eta \Delta y_\xi$. Note that $\mathbbm{1}_{\omega} h$ is localized in $\omega$ because of the support properties of $\eta$. 

Therefore, $y_0$ vanishes in $(0,T) \times \Theta$, hence the associated function $z_0$ such that $(y_0,z_0)$ satisfies \eqref{y-Omega-0-h} is identically zero. In particular, according to \eqref{Diff-J-V-2}, we immediately have the exact insensitizing property \eqref{Exact-Insentitizing}.

%%%%%%%%%%%%%%%%%%%%%%%%%%%%
\subsection{Proof of Theorem \ref{Thm-insensitizing-Boundary}}
\label{Sec-Proof-Thm-insensitizing-Boundary}
%%%%%%%%%%%%%%%%%%%%%%%%%%%%

We start by introducing open sets $\omega_0$, $\omega_1$, $\omega_2$, and $\omega_3$ such that
\begin{equation*}
	\partial \Theta \subset \omega_0 \Subset  \omega_1 \Subset  \omega_2 \Subset  \omega_3 \Subset  \omega,
	% \qquad \hbox{ with } \omega_3 \cap \Theta \subset \Theta.
\end{equation*}
which is possible thanks to Assumption \eqref{Assumption-insensitizing-Boundary}.
We also introduce a smooth function $\eta_{23} = \eta_{23}(x)$  taking value $1$ in $\Omega \setminus \overline{\omega_3}$ and vanishing in $\overline{\omega_2}$.
 
We fix $\xi \in L^2(0,T; L^2(\R^d))$, and we introduce the solution $y_\xi$ to
\begin{equation*}
%	\label{Heat-Eq-y-f}
	\left\{
		\begin{array}{ll}
			\partial_t y_\xi - \Delta y_\xi = \eta_{23} \xi, \quad & (t,x) \in (0,T) \times \Omega_0, 
			\\
			y_\xi(t, x) = 0, \quad & (t,x) \in (0,T) \times \partial \Omega_0, 
			\\
			y_\xi(0,x) = 0, \quad & x \in \Omega, 
		\end{array}
	\right.
\end{equation*}
and the solution  $z_\xi$ to
\begin{equation}
	\label{Heat-Eq-Adjoint-z-f}
	\left\{
		\begin{array}{ll}
			- \partial_t z_\xi - \Delta z_\xi   = \eta_{12} y_\xi  \mathbbm{1}_{\Theta} , \quad & (t,x) \in (0,T) \times \Omega_0, 
			\\
			z_\xi(t, x) = 0, \quad & (t,x) \in (0,T) \times \partial \Omega_0, 
			\\
			z_\xi(T,x) = 0, \quad & x \in \Omega_0. 
		\end{array}
	\right.
\end{equation}
where $\eta_{12} = \eta_{12}(x)$ is a smooth function taking value $1$ in $\Omega_0 \setminus \overline{\omega_2}$ and vanishing in $\overline{\omega_1}$.

Then, we  introduce a smooth function $\eta_{01} = \eta_{01}(x)$ such that $\eta_{01}$ vanishes in $ \omega_0$ and equal to $1$ in $\Omega_0 \setminus \omega_1$, so that $\mathbbm{1}_\Theta \eta_{01}$ is actually a smooth function taking value $1$ in $\Theta \setminus \omega_1$ and vanishing in $\omega_0 \cup (\Omega_0 \setminus \Theta)$. Then $ z_0(t,x)= \mathbbm{1}_\Theta(x) \eta_{01}(x) z_\xi(t,x)$ satisfies
\begin{equation*}
%	\label{Heat-Eq-Adjoint-q-final}
	\left\{
		\begin{array}{ll}
			- \partial_t z_0 - \Delta z_0   = \eta_{01} \eta_{12} y_\xi  \mathbbm{1}_{\Theta} 
								- [\Delta, \mathbbm{1}_\Theta \eta_{01}] z_\xi
								, \quad & (t,x) \in (0,T) \times \Omega_0, 
			\\
			z_0(t, x) = 0, \quad & (t,x) \in (0,T) \times \partial \Omega_0, 
			\\
			z_0(T,x) = 0, \quad & x \in \Omega_0, 
		\end{array}
	\right.
\end{equation*}
and, by construction, 
\[
	z_0(t,x) = 0 \hbox{ for all } (t,x) \in (0,T) \times (\Omega_0 \setminus \Theta).
\]
Now, we remark that by construction, $[\Delta, \mathbbm{1}_\Theta \eta_{01}] z_\xi$ is localized in $\omega_1 \cap \Theta$. Besides, local regularity results for \eqref{Heat-Eq-Adjoint-z-f} imply that $z_\xi$ is $\mathscr{C}^2([0,T] \times \overline{\omega_1 \cap \Theta})$. We then take 
\[
	y_0(t,x) = \eta_{01}(x) \eta_{12}(x) y_\xi(t,x) - [\Delta, \mathbbm{1}_\Theta \eta_{01} ] z_\xi(t,x), \quad \hbox{ for all } (t,x) \in (0,T) \times \Omega_0, 
\]
which satisfies \eqref{y-Omega-0-h} for $h$ given for $(t,x) \in (0,T) \times \Omega_0$ by
\[
	h(t,x) = (\eta_{01}(x) \eta_{12}(x) \eta_{23}(x)-1) \xi(t,x) - [\Delta, \eta_{01} \eta_{12}] y_\xi(t,x) - (\partial_t - \Delta) ([\Delta, \mathbbm{1}_\Theta \eta_{01}] z_\xi)(t,x). 
\]	
This control function $h$ is localized in $(0,T) \times \omega$ due to the conditions on the support of $\eta_{01}$, $\eta_{12}$, $\eta_{23}$. This concludes the proof of Theorem \ref{Thm-insensitizing-Boundary}.

\subsection{Proof of Theorem \ref{Thm-Negative-Theta-Omega}}
\label{Sec-Proof-Thm-Negative-Theta-Omega}

Firstly, we consider a function $g \in \mathscr{C}^4(\partial\Omega_0)$ such that $g$ is nowhere $\mathscr{C}^5(\partial \Omega_0)$ (such functions form a dense set in the sense of Baire of $\mathscr{C}^4(\partial \Omega_0)$ and thus exist). 

Then, we introduce a function $q_* \in H^4(\Omega_0)$ such that 
\begin{align*}
	&
	q_* (x) = 0\ \hbox{ for } x \in \partial \Omega_0, 
	\qquad
	&&\partial_n q_* (x) = g(x)\ \hbox{ for } x \in \partial \Omega_0, 
	\\
	&
	\Delta q_* (x) = 0\ \hbox{ for } x \in \partial \Omega_0, 
	\qquad
	&&\partial_n \Delta q_* (x) = g(x)\ \hbox{ for } x \in \partial \Omega_0, 
\end{align*}
{whose existence is guaranteed by classical trace theorems, see e.g. \cite[Chap.1, Theorem 8.3]{LionsMagenes},} and we choose a smooth non-negative function $\eta = \eta(t)$ such that $\eta(0) = \eta(T) =\eta'(0) = 0$ and $\int_0^T \eta(t)^2 \, dt = 1$. 

Then, we set, for $(t,x) \in (0,T) \times \Omega_0$, 
\begin{align*}
	& z_\xi(t,x) = \eta(t) q_*(x), 
	\quad
	y_\xi(t,x) = {- \partial_t z_\xi(t,x) - \Delta z_\xi(t,x)}
		=  - \eta'(t) q_*(x) - \eta(t) \Delta q_*(x), 
	\\
	& \xi(t,x) = \partial_t y_\xi(t,x) - \Delta y_\xi(t,x). 
\end{align*}
{Note that $\xi \in L^2(0,T; L^2( \Omega_0))$ since $q_* \in H^4(\Omega_0)$.}

Assume that we can solve the exact insensitizing problem for this choice of $\xi$. Hence, from \eqref{Exact-Insentitizing} and formula \eqref{Diff-J-V-2}, there exists $h \in L^2(0,T; L^2(\omega))$ such that 
\begin{equation}
	\label{Equivalent-y-OmegaV-condition-Gamma}
	\forall x \in \partial \Omega_0, \quad \int_0^T \partial_n y_0(t,x) \partial_n z_0(t,x) \, dt = 0, 
\end{equation}
where $(y_0, z_0)$ solves \eqref{y-Omega-0-h}.
Now, we decompose $y_0$ as
\begin{equation*}
	y(t,x) = y_\xi (t,x) + y_h(t,x), 
\end{equation*}
where $y_h$ is the solution of \eqref{cou1bis}, and $z_0$ as 
\begin{equation*}
	z_0(t,x) = z_\xi(t,x) + z_h(t,x),  \quad (t,x) \in (0,T) \times \Omega_0,
\end{equation*}
where $z_h$ solves \eqref{cou2bis}.

From \eqref{Equivalent-y-OmegaV-condition-Gamma},  for all $x \in \partial \Omega_0$,
\begin{align*}
	0 = & \int_0^T (\partial_n y_\xi + \partial_n y_h) (\partial_n z_\xi + \partial_n z_h) \, dt
	\\
	= &
	- (g(x))^2 + g(x) \int_0^T \left(\eta(t) \partial_n y_h(t,x)  - (\eta' (t)+ \eta(t)) \partial_n z_h(t,x) \right) \, dt
	+ \int_0^T \partial_n y_h(t,x)\, \partial_n z_h(t,x) \, dt, 
\end{align*}
where we used that 
\[
	\partial_n y_\xi = - ( \eta + \eta') g, \quad \partial_n z_\xi = \eta g. 
\]
Since $\omega \Subset  \Omega_0$, {the regularizing properties of the heat equation imply that $y_h$ is smooth close to the boundary $[0,T] \times \partial \Omega_0$, and thus so is $z_h$. Therefore, the quantities} 
\[
	a_0(x) = \int_0^T \left(\eta(t) \partial_n y_h(t,x)  - (\eta' (t)+ \eta(t)) \partial_n z_h(t,x) \right) \, dt
	\quad \hbox{ and } \quad 
	a_1(x) =   \int_0^T \partial_n y_h(t,x)\, \partial_n z_h(t,x) \, dt, 
\]
are smooth ($\mathscr{C}^\infty$) in $\partial \Omega_0$.  Since for all $x \in \partial \Omega_0$, $g(x)$ is a real root to the polynomial $-X^2 + X a_0(x) + a_1(x)$, we necessarily have that for all $x \in \partial \Omega_0$, $a_0(x)^2 + 4 a_1(x) \geq 0$ and for all $x \in \partial \Omega_0$,
\[
	g(x) \in \left\{
		 \frac{1}{2} \Big(a_0(x) + \sqrt{a_0(x)^2 + 4 a_1(x)}\Big), \, 
		 \frac{1}{2} \Big(a_0(x) - \sqrt{a_0(x)^2 + 4 a_1(x)}\Big)
		 \right\}.
\]
In particular, if there exists $x_0 \in \partial \Omega$ such that $a_0(x_0)^2 > 4 a_1(x_0)$, since $g$ is continuous, there is a sign $s \in \{-1,1\}$ such that in a neighbourhood of $x_0$ (in $\partial \Omega$) in which $a_0^2 + 4 a_1$ stays positive, 
\[
	g(x) =  \frac{1}{2} \Big(- a_0(x) + s \sqrt{a_0(x)^2 {+}4 a_1(x)}\Big),
\]
implying in particular that $g$ is smooth ($\mathscr{C}^\infty$) in a neighbourhood of $x_0$, which contradicts the choice of $g$. 

Thus, for all $x \in \partial \Omega_0$, we should have $a_0(x)^2 + 4 a_1(x) = 0 $, so that $g(x)  = -a_0(x)/2$. But this would again imply that $g$ is smooth in $\partial \Omega_0$, thus contradicting the choice of $g$. 

We have thus obtained a contradiction. There cannot be any control $h \in L^2(0,T; L^2(\omega))$ such that the condition \eqref{Equivalent-y-OmegaV-condition-Gamma} holds.

%\appendix
%\label{app1}
%Our goal is to prove the following Lemma.
%\begin{lemma}
%Consider a weak solution of 

%%%%%%%%%%%%%%%%%%%%%
%%%%%%%%%%%%%%%%%%%%%%%%%%%%%%%%%%%%%%%%%%%%%%%%%%%%%%%%%%%%%%%%%%%%%%%%%%%%%%%%%%%%%%%%%%%%%%%%%%%%%%%%%%%%%%%%%%%%%%%%%%%%%%%%%%%%%%%%%%%%%%%%%%%%%%%%%%%%%%%%%%%%%%%%%%%%%%%%%%%%%%%%%%%%%%%%%%%%%%%%%%%%%%%%%%%%%%%%%%%%%%%%%%%%%%%%%%%%%%%%%%%%%%%%%%%%%%%%%%%%%%%%%%%%%%%%%%%%%%%%%%%%%%%%%%%%%%%%%%%%%%%%%%%%%%%%%%%%%%%%%%%%%%%%%%%%%%%%%%%%%%%%%%%%%%%%%%%%%%%%%%%%%%%%%%%%%%%%%%%%%%%%%%%%%%%%%%%%%%%

%%%%%%%%%%%%%%%%%%%%%%%%%%%%%%%%%%%%%%%%%%%
%
\bibliographystyle{plain}
%\bibliography{/Users/servedoza/Desktop/Biblio} 

\begin{thebibliography}{}

\end{thebibliography}


\begin{thebibliography}{10}

\bibitem{Alabau-Boussouira}
F.~Alabau-Boussouira.
\newblock Insensitizing exact controls for the scalar wave equation and exact controllability of 2-coupled cascade systems of PDE's by a single control.
\newblock {\em Math. Control Signals Systems}, 26, no. 1, 1--46, 2014.
\bibitem{mousquetaires}F.~Ammar Khodja, A.~ Benabdallah, M.~ Gonz\'alez-Burgos and L. de Teresa. 
\newblock
New phenomena for the null controllability of parabolic systems: minimal time and geometrical dependence. (English summary)
\newblock{\em J. Math. Anal. Appl.}, 444, no. 2, 1071--1113, 2016.

\bibitem{Beauchard-Marbach-2017}
K.~Beauchard and F.~Marbach.
\newblock Quadratic obstructions to small-time local controllability for
  scalar-input systems.
\newblock {\em J. Differential Equations}, 264(5):3704--3774, 2018.

\bibitem{Bodart-Fabre}
O.~Bodart and C.~Fabre.
\newblock Controls insensitizing the norm of the solution of a semilinear heat
  equation.
\newblock {\em J. Math. Anal. Appl.}, 195(3):658--683, 1995.

\bibitem{Bodart-GB}
O.~Bodart, M.~Gonz\'alez-Burgos and R.~P\'erez-Garc\'ia.
\newblock Insensitizing controls for a heat equation with a nonlinear term involving the state and the gradient.
\newblock {\em Nonlinear Anal.}, 57, no. 5-6, 687--711, 2004.

\bibitem{Boyer-HS-dT}
F.~Boyer., V. Hern\'andez-Santamar\'a and L.~de~Teresa
\newblock Insensitizing controls for a semilinear parabolic equation: a numerical approach.
\newblock {\em Math. Control Relat. Fields}, 9, no. 1, 117–158, 2019.

\bibitem{Brezis}
H.~Brezis.
\newblock {\em Functional analysis, {S}obolev spaces and partial differential
  equations}.
\newblock Springer, 2011.

\bibitem{Carreno}
N.~Carre\~no.
\newblock Insensitizing controls for the Boussinesq system with no control on the temperature equation.
\newblock {\em Adv. Differential Equations.}, 22, no. 3-4, 235--258, 2017.

\bibitem{Carreno-Guer-Gueye}
N.~Carre\~no, S.~Guerrero and M.~Gueye.
\newblock Insensitizing controls with two vanishing components for the three-dimensional Boussinesq system.
\newblock {\em ESAIM Control Optim. Calc. Var.}, 21, no. 1, 73--100, 2015.

\bibitem{Carreno-Gueye}
N.~Carre\~no and M.~Gueye.
\newblock Insensitizing controls with one vanishing component for the Navier-Stokes system.
\newblock {\em J. Math. Pures Appl. (9)}, 101, no. 1, 27–53, 2014.

\bibitem{Chowdhury-Erv-2019}
S.~Chowdhury and S.~Ervedoza.
\newblock Open loop stabilization of incompressible {N}avier-{S}tokes equations
  in a 2d channel using power series expansion.
\newblock {\em J. Math. Pures Appl. (9)}, 130:301--346, 2019.

\bibitem{CNL}
 J.-M. Coron.
 \newblock Control and nonlinearity. 
 \newblock{\em Mathematical Surveys and Monographs}, 136. American Mathematical Society, 2007.
 
\bibitem{Delfour-Zolesio}
M.~C. Delfour and J.-P. Zol\'{e}sio.
\newblock {\em Shapes and geometries}, volume~22 of {\em Advances in Design and
  Control}.
\newblock Society for Industrial and Applied Mathematics (SIAM), Philadelphia,
  PA, second edition, 2011.
\newblock Metrics, analysis, differential calculus, and optimization.

\bibitem{Ervedoza-Note-2020}
S.~Ervedoza.
\newblock Control issues and linear projection constraints on the control and
  on the controlled trajectory.
\newblock {\em North-West. Eur. J. Math.}, 6:165--198, 2020.

\bibitem{Evans}
L. C. Evans.
\newblock Partial differential equations. 
\newblock{\em Graduate Studies in Mathematics}, 19. American Mathematical Society, Providence, RI, 1998.

\bibitem{FernandezCaraZuazua1}
E.~Fern{\'a}ndez-Cara and E.~Zuazua.
\newblock The cost of approximate controllability for heat equations: the
  linear case.
\newblock {\em Adv. Differential Equations}, 5(4-6):465--514, 2000.

\bibitem{FursikovImanuvilov}
A.~V. Fursikov and O.~Y. Imanuvilov.
\newblock {\em Controllability of evolution equations}, volume~34 of {\em
  Lecture Notes Series}.
\newblock Seoul National University Research Institute of Mathematics Global
  Analysis Research Center, Seoul, 1996.

\bibitem{Guerrero}
S.~Guerrero.
\newblock Controllability of systems of Stokes equations with one control force: existence of insensitizing controls.
\newblock {\em Ann. Inst. H. Poincaré Anal. Non Linéaire} 24, no. 6, 1029--1054, 2007.

\bibitem{Gueye}
M.~Gueye.
\newblock Insensitizing controls for the Navier-Stokes equations.
\newblock {\em Ann. Inst. H. Poincaré Anal. Non Linéaire} 30, no. 5, 825--844, 2013.

\bibitem{Henrot-Pierre-2005}
A.~Henrot and M.~Pierre.
\newblock {\em Shape Variation and Optimization}, volume~28 of {\em Tracts in
  Mathematics}.
\newblock European Mathematical Society, Z\"urich, 2018.

\bibitem{Hormander-2003}
L.~H\"ormander.
\newblock {\em The analysis of linear partial differential operators. I.
Distribution theory and Fourier analysis.} Reprint of the second (1990) edition. Classics in Mathematics. Springer-Verlag, Berlin, 2003. 

\bibitem{Kavian-De-Teresa}
O.~Kavian and L.~de~Teresa.
\newblock Unique continuation principle for systems of parabolic equations.
\newblock {\em ESAIM Control Optim. Calc. Var.}, 16(2):247--274, 2010.

\bibitem{Lions-1990}
J.-L. Lions.
\newblock Quelques notions dans l'analyse et le contr\^ole de syst\`emes \`a donn\'ees incompl\`etes. (French) [Some notions in the analysis and control of systems with incomplete data] 
\newblock Proceedings of the XIth Congress on Differential Equations and Applications/First Congress on Applied Mathematics (Spanish) (M\'alaga, 1989), 43–54, Univ. M\'alaga, M\'alaga, 1990. 

\bibitem{Lions-1992}
J.-L. Lions.
\newblock Remarks on approximate controllability.
\newblock {\em J. Anal. Math.}, 59:103--116, 1992.
\newblock Festschrift on the occasion of the 70th birthday of Shmuel Agmon.

\bibitem{Lions-Sentinelles-1992}
J.-L. Lions.
\newblock {\em Sentinelles pour les syst\`emes distribu\'{e}s \`a donn\'{e}es
  incompl\`etes}, volume~21 of {\em Recherches en Math\'{e}matiques
  Appliqu\'{e}es [Research in Applied Mathematics]}.
\newblock Masson, Paris, 1992.

\bibitem{LionsMagenes}
J.-L. Lions and E.~Magenes.
\newblock {\em Probl\`emes aux limites non homog\`enes et applications. {V}ol.
  1}.
\newblock Travaux et Recherches Math\'ematiques, No. 17. Dunod, Paris, 1968.

\bibitem{LionsMagenes2}
J.-L. Lions and E.~Magenes.
\newblock {\em Probl\`emes aux limites non homog\`enes et applications. {V}ol.
  2}.
\newblock Travaux et Recherches Math\'ematiques, No. 17. Dunod, Paris, 1968.

\bibitem{Lissy-Privat-Simpore}
P.~Lissy, Y.~Privat, and Y.~Simpor\'{e}.
\newblock Insensitizing control for linear and semi-linear heat equations with
  partially unknown domain.
\newblock {\em ESAIM Control Optim. Calc. Var.}, 25:Art. 50, 21, 2019.

%\bibitem{Mizo}
%S.~Mizohata.
%\newblock Unicit\'e du prolongement des solutions pour quelques op\'erateurs diff\'erentiels paraboliques.
%\newblock{\em Mem. Coll. Sci. Univ. Kyoto Ser. A. Math.}, 31, 219–239, 1958.

\bibitem{Micu-Ortega-deTeresa}
S.~Micu, J.~H. Ortega, and L.~de~Teresa.
\newblock An example of {$\epsilon$}-insensitizing controls for the heat
  equation with no intersecting observation and control regions.
\newblock {\em Appl. Math. Lett.}, 17(8):927--932, 2004.

\bibitem{Tebou1}
L.~Tebou.
\newblock Some results on the controllability of coupled semilinear wave equations: the desensitizing control case.
\newblock {\em SIAM J. Control Optim.},  49 (2011), no. 3, 1221--1238.

\bibitem{Tebou2}
L.~Tebou.
\newblock Locally distributed desensitizing controls for the wave equation.
\newblock {\em C. R. Math. Acad. Sci. Paris},  346 (2008), no. 7-8, 407--412.

\bibitem{deTeresa-00}
L.~de Teresa.
\newblock Insensitizing controls for a semilinear heat equation.
\newblock {\em Comm. Partial Differential Equations} 25, no. 1-2, 39--72, 2000.

\bibitem{dT-Zuazua-09}
L.~de Teresa and E.~Zuazua.
\newblock Identification of the class of initial data for the insensitizing control of the heat equation.
\newblock {\em Commun. Pure Appl. Anal. (8)}, no. 1, 457--471, 2009.

\red{
\bibitem{TWBook}
M.~Tucsnak and G.~Weiss.
\newblock {\em Observation and control for operator semigroups}.
\newblock Birkh\"auser Advanced Texts: Basler Lehrb\"ucher. [Birkh\"auser
  Advanced Texts: Basel Textbooks]. Birkh\"auser Verlag, Basel, 2009.
}

\bibitem{Zuazua-97}
E.~Zuazua.
\newblock Finite-dimensional null controllability for the semilinear heat
  equation.
\newblock {\em J. Math. Pures Appl. (9)}, 76(3):237--264, 1997.

\end{thebibliography}

%
%
\end{document}